# WEAK SOLUTIONS FOR FORWARD–BACKWARD SDES—A MARTINGALE PROBLEM APPROACH

By Jin Ma,[1] Jianfeng Zhang[2] and Ziyu Zheng[3]

*University of Southern California and Barclays Capital*

In this paper, we propose a new notion of *Forward–Backward Martingale Problem* (FBMP), and study its relationship with the weak solution to the forward–backward stochastic differential equations (FBSDEs). The FBMP extends the idea of the well-known (forward) martingale problem of Stroock and Varadhan, but it is structured specifically to fit the nature of an FBSDE. We first prove a general sufficient condition for the existence of the solution to the FBMP. In the Markovian case with uniformly continuous coefficients, we show that the weak solution to the FBSDE (or equivalently, the solution to the FBMP) does exist. Moreover, we prove that the uniqueness of the FBMP (whence the uniqueness of the weak solution) is determined by the uniqueness of the viscosity solution of the corresponding quasilinear PDE.

**1. Introduction.** The theory of backward stochastic differential equations (BSDEs for short) has matured tremendously since the seminal work of Pardoux and Peng [24]. The fundamental well-posedness of BSDEs with various conditions on the coefficients as well as terminal conditions have been studied extensively, which can be found in a large amount of literature. A commonly used list of reference include the books of El Karoui and Mazliak [9] and Ma and Yong [18] for the basic theory of BSDEs, and the survey paper of El Karoui, Peng and Quenez [10] for the applications of BSDEs to mathematical finance. It is worth noting that almost all the existing works on BSDEs or its extension Forward–Backward SDEs (FBSDEs) are exclusively considered in the realm of "strong solutions," and a missing

Received February 2006; revised June 2007.
[1]Supported in part by NSF Grants 05-05427 and 08-06017.
[2]Supported in part by NSF Grant 04-03575 and 06-31366.
[3]Supported in part by NSF Grant 03-06233.

piece of puzzle in the theory of BSDEs seems to have been the notion of "weak-solution." Such a notion, although extremely conceivable and tempting from both theoretical and practical point of views, has not been fully explored.

In a recent paper, Antonelli and Ma [1] introduced the notion of weak solution to a class of FBSDEs. In that paper, it was shown that some standard results regarding the relations among weak solution, strong solution and different types of uniqueness still hold. However, the results in that paper were a far cry from a systematic study for weak solutions. In particular, the authors were not able to address the core issue regarding the uniqueness. Similar topics were studied later by Buckdahn, Engelbert and Rascanu [4], with a more general definition and more extended investigation. But the issue of uniqueness remains. Independent of our work, recently Delarue and Guatteri [8] established the existence and uniqueness of weak solutions for a class of Markovian FBSDEs by using elegantly the decoupling strategy in the Four Step Scheme (cf. [16]). However, they require the coefficients to be Lipschitz continuous in the backward components, and thus roughly speaking their FBSDE is weak only in the forward component. To our best knowledge, so far there has not been any work trying to address the issue of uniqueness in law for a true BSDE/FBSDE; and it is our hope that this paper could be the first step in that direction.

Our first goal of this paper is to find an appropriate definition of a "backward" version of the *martingale problem* associated to the weak solution. We shall follow naturally the idea of the forward martingale problem initiated by Stroock and Varadhan (cf., e.g., [26]), and recast the FBSDE in terms of some fundamental martingales, which then leads to the notion of the *Forward–Backward Martingale Problem* (FBMP). Such a notion extends the usual martingale problem and it is equivalent to the weak solution defined in [1]. Our objective then is to prove the existence and uniqueness of the solution to FBMP, whence those of weak solution. Given the large amount of recent studies on the existence of (strong) solutions to BSDEs/FBSDEs with less-regular coefficients, notably the works of [3, 6, 14] and [15], to mention a few, we are particularly interested in finding a unified method that works for high dimensional FBSDEs with nonsmooth coefficients for which a strong solution is less likely to exist. We shall first prove a general sufficient condition for the existence of solution to FBMP, using mainly some weak convergence arguments under Meyer–Zheng topology. A key element in the sufficient condition could essentially be understood as a certain type of tightness criterion for processes with paths in an $L^2$ space, which shall be further explored in our future publications. We will then show that such sufficient condition can be verified in a Markovian case assuming that all the coefficients are bounded and uniformly continuous.

The last part of this paper deals with the main issue: the uniqueness of the solution to the FBMP. We note that to date the main difficulties in the discussion has always been the martingale integrand in the BSDE (the process $Z$), because in general one does not have a workable canonical space for this process. In fact, although in many cases the process $Z$ is càdlàg or even continuous (see, e.g., [20]), such path regularity is by no means clear a priori. However, it is noted that if all the coefficients are Hölder continuous, one can show that the martingale integrand can be treated as a function of the forward components of the solution, owing to the idea of the Four Step Scheme of [16]. This fact, together with the procedure we used to prove the existence, shows that at least one weak solution can be built using only the path spaces of the continuous components of the solution. This result becomes more significant when we establish the uniqueness, since it essentially eliminated the subtlety caused by the canonical spaces. Our uniqueness proof is originated from the idea of "*method of optimal control*" for solving an FBSDE (see [17, 18]). Although it may not be intuitive due to the technicalities involved, the basic idea is to investigate a variation of the notion of "nodal set" in [17], so as to show that uniqueness of the viscosity solutions to the corresponding quasilinear PDE implies the uniqueness of the solution to the FBMP (whence the weak solution). We should note that in this paper we are still not able to prove the uniqueness in the most general sense, but we believe that our method has a potential to be applied to more general FBSDEs, and the uniqueness should hold in a much wider class of weak solutions. We hope to be able to address the issue in our future publications.

The rest of the paper is organized as follows. In Section 2, we give the preliminaries, recall the definition of weak solution and introduce the notion of an FBMP. In Section 3, we prove the general sufficient condition for the existence of the solution to FBMP. In Section 4, we consider the Markovian case. Finally, in Section 5, we prove the uniqueness of the solution to FBMP.

**2. Preliminaries.** In this section, we give the basic probabilistic set up, recall the definition of weak solution of an FBSDE and introduce the notion of Forward–Backward Martingale Problems (FBMPs).

For any Euclidean space $\mathbb{R}^k$, regardless of its dimension, we denote its norm by $|\cdot|$. We denote $\mathbb{C}([0,T];\mathbb{R}^k)$ to be the space of all $\mathbb{R}^k$-valued continuous functions endowed with the sup-norm; and $\mathbb{D}([0,T];\mathbb{R}^k)$ to be the space of all $\mathbb{E}$-valued càdlàg functions endowed with the Skorohod topology (see, e.g., [11]). When $k=1$, we may omit $\mathbb{R}$ in the notation.

For a given finite time horizon $[0,T]$, we say that a quintuple $(\Omega, \mathcal{F}, P, \mathbf{F}, W)$ is a *standard set-up* if $(\Omega, \mathcal{F}, P)$ is a complete probability space; $\mathbf{F} \stackrel{\triangle}{=} \{\mathcal{F}_t\}_{t\in[0,T]}$ is a filtration satisfying the *usual hypotheses* (see, e.g., [25]); and $W$ is an $\{\mathcal{F}_t\}$-Brownian motion. In particular, if $\mathcal{F}_t = \mathcal{F}_t^W$, the natural filtration generated by the Brownian motion $W$, augmented by all the $P$-null sets of $\mathcal{F}$ and satisfying the usual hypotheses, then we say that the standard set-up is *Brownian*.

A. WEAK SOLUTION OF FBSDEs. Let us consider the following forward–backward SDE:

$$
(2.1) \quad \begin{cases} X_t = x + \int_0^t b(s, (X)_s, Y_s, Z_s)\, ds + \int_0^t \sigma(s, (X)_s, Y_s, Z_s)\, dW_s, \\ Y_t = g((X)_T) + \int_t^T h(s, (X)_s, Y_s, Z_s)\, ds - \int_t^T Z_s\, dW_s. \end{cases}
$$

Here, $(X_t, Y_t, Z_t, W_t) \in \mathbb{R}^n \times \mathbb{R}^m \times \mathbb{R}^{m \times d} \times \mathbb{R}^d$, and the functions $b$, $h$, $\sigma$ and $g$ are functions with appropriate dimensions. We note, in particular, that the coefficient $b$ is *a progressively measurable function* defined on $[0, T] \times \mathbb{C}([0, T], \mathbb{R}^n) \times \mathbb{R}^m \times \mathbb{R}^{m \times d}$ with valued in $\mathbb{R}^n$, and $(X)_t$ denotes the path of $X$ up to time $t$. More precisely, for each $t \in [0, T]$ and $(y, z) \in \mathbb{R}^m \times \mathbb{R}^{m \times d}$, the mapping $\mathbf{x} \mapsto b(t, (\mathbf{x})_t, y, z)$ is measurable with respect to the $\sigma$-field $\mathcal{B}_t(\mathbb{C}([0, T]; \mathbb{R}^n))$, where $\mathcal{B}_t(\mathbb{C}([0, T]; \mathbb{R}^n)) \stackrel{\triangle}{=} \sigma\{\mathbf{x}(t \wedge \cdot) : \mathbf{x} \in \mathbb{C}([0, T]; \mathbb{R}^n)\}$ (cf., e.g., [13]). The coefficients $\sigma$, $h$ and $g$ should be understood in a similar way. It is known that (cf., e.g., [18]) an adapted (strong) solution to the FBSDE (2.1) is usually understood as a triplet of processes $(X, Y, Z)$ defined on any given Brownian set-up such that (2.1) holds $P$-almost surely. The following definition of *weak solution* is proposed in [1].

DEFINITION 2.1. A standard set-up $(\Omega, \mathcal{F}, P, \{\mathcal{F}_t\}, W)$ along with a triplet of processes $(X, Y, Z)$ defined on this set-up is called a weak solution of (2.1) if:

(i) the processes $X, Y$ are continuous, and all processes $X$, $Y$, $Z$ are $\mathcal{F}_t$-adapted;

(ii) denoting $f_t = f(t, (X)_t, Y_t, Z_t)$ for $f = b, \sigma, h$, it holds that

$$P\left\{\int_0^T (|b_t| + |\sigma_t|^2 + |h_t| + |Z_t|^2)\, ds + |g((X)_T)| < \infty\right\} = 1.$$

(iii) $(X, Y, Z)$ verifies (2.1) $P$-a.s.

We remark here that unlike a "strong solution," a weak solution relaxed the most fundamental requirement for a BSDE, that is, the set-up be Brownian. But instead, it requires the flexibility of the set-up for each solution, similar to the forward SDE case. We should point out that in [1] it is shown that the weak solution of FBSDE (2.1) exists under very mild conditions, and that there does exist a weak solution that is not "strong."

REMARK 2.2. Although in the basic setting of FBSDE (2.1), the coefficients are seemingly "deterministic," it can be easily extended to the "random coefficients" case. For instance, if we add the canonical Brownian motion $W$ into the equation, and consider $(W, X)$ as the forward component, then we can allow the coefficients to have the form

$$(2.2) \quad f(t, \omega, (X)_t, Y_t Z_t) \stackrel{\triangle}{=} f(t, (W)_t, (X)_t, Y_t, Z_t), \qquad f = b, \sigma, h,$$

and the FBSDE (2.1) has nonanticipating random coefficients. In fact, our general existence result Theorem 3.1 holds true for general FBSDEs with coefficients of the form (2.2). However, at this stage, we feel that it is more convenient to consider (2.1) in the given form so as to avoid further complication in the proof of the uniqueness. We should note that even in the standard (forward) martingale problem (cf. [26]), the component $W$ is not involved directly. □

B. FORWARD–BACKWARD MARTINGALE PROBLEM. Before we define the martingale problem, let us give a detailed description of a "canonical set-up" on which our discussion will be carried out. Define

$$(2.3) \quad \Omega^1 \stackrel{\triangle}{=} \mathbb{C}([0,T]; \mathbb{R}^n); \qquad \Omega^2 \stackrel{\triangle}{=} \mathbb{C}([0,T]; \mathbb{R}^m); \qquad \Omega \stackrel{\triangle}{=} \Omega^1 \times \Omega^2,$$

where $\Omega^1$ denotes the path space of the forward component $X$ and $\Omega^2$ the path space of the backward component $Y$ of the FBSDE, respectively.

Next, we define the canonical filtration by $\mathcal{F}_t \stackrel{\triangle}{=} \mathcal{F}_t^1 \otimes \mathcal{F}_t^2$, $0 \le t \le T$, where $\mathcal{F}_t^i \stackrel{\triangle}{=} \sigma\{\omega^i(r \wedge t) : r \ge 0\}$, $i = 1, 2$. We denote $\mathcal{F} \stackrel{\triangle}{=} \mathcal{F}_T$ and $\mathbf{F} \stackrel{\triangle}{=} \{\mathcal{F}_t\}_{0 \le t \le T}$.

In what follows, we denote the generic element of $\Omega$ by $\omega = (\omega^1, \omega^2)$, and denote the canonical processes on $(\Omega, \mathcal{F})$ by

$$\mathbf{x}_t(\omega) \stackrel{\triangle}{=} \omega^1(t) \quad \text{and} \quad \mathbf{y}_t(\omega) \stackrel{\triangle}{=} \omega^2(t), \quad t \ge 0.$$

Finally, let $\mathcal{P}(\Omega)$ be all the probability measures defined on $(\Omega, \mathcal{F})$, endowed with the Prohorov metric.

To simplify presentation, we first assume that $\sigma = \sigma(t, (\mathbf{x})_t, y)$. Here, we abuse the notation $\mathbf{x}$ by denoting elements of $\mathbb{C}([0,T), \mathbb{R}^n)$ instead of the canonical process. [The case $\sigma = \sigma(t, (\mathbf{x})_t, y, z)$ is a little more complicated; we address it in Remark 2.4 below.] Further, for $f = b, h$, we denote $\hat{f}(t, (\mathbf{x})_t, y, z) = f(t, (\mathbf{x})_t, y, z\sigma(t, (\mathbf{x})_t, y))$, and let $a = \sigma\sigma^T$. We give the following definition for a *forward–backward martingale problem*.

DEFINITION 2.3. Let $b$, $\sigma$, $h$ and $g$ be given. For any $x \in \mathbb{R}^n$, a solution to the forward–backward martingale problem with coefficients $(b, \sigma, h, g)$ [FBMP$_{x,T}(b, \sigma, h, g)$ for short] is a pair $(\mathbb{P}, \mathbf{z})$, where $\mathbb{P} \in \mathcal{P}(\Omega)$, and $\mathbf{z}$ is a $\mathbb{R}^{m \times n}$-valued predictable process defined on the filtered canonical space $(\Omega, \mathcal{F}, \mathbf{F})$, such that following properties hold:

(i) the processes

(2.4)
$$M_{\mathbf{x}}(t) \triangleq \mathbf{x}_t - \int_0^t \widehat{b}(r, (\mathbf{x})_r, \mathbf{y}_r, \mathbf{z}_r)\, dr \quad \text{and}$$

$$M_{\mathbf{y}}(t) \triangleq \mathbf{y}_t + \int_0^t \widehat{h}(r, (\mathbf{x})_r, \mathbf{y}_r, \mathbf{z}_r)\, dr$$

are both $(\mathbb{P}, \mathbf{F})$-martingales for $t \in [0, T]$;

(ii) $[M_{\mathbf{x}}^i, M_{\mathbf{x}}^j](t) = \int_0^t a_{ij}(r, (\mathbf{x})_r, \mathbf{y}_r)\, dr$, $t \in [0, T]$, $i, j = 1, \ldots, n$;
(iii) $M_{\mathbf{y}}(t) = \int_0^t \mathbf{z}_r\, dM_{\mathbf{x}}(r)$, $t \in [0, T]$.
(iv) $\mathbb{P}\{\mathbf{x}_0 = x\} = 1$ and $\mathbb{P}\{\mathbf{y}_T = g((\mathbf{x})_T)\} = 1$.

We note that by (iii) we imply that the quadratic variation of $M_{\mathbf{y}}$ is absolutely continuous with respect to the quadratic variation of $M_{\mathbf{x}}$, thus in the definition we require implicitly

$$\mathbb{P}\bigg\{\int_0^T |\mathbf{z}_t a(t, (\mathbf{x})_t, \mathbf{y}_t) \mathbf{z}_t^T|_{\mathbb{R}^{m \times m}}^2\, dt < \infty\bigg\} = 1.$$

REMARK 2.4. The case when $\sigma = \sigma(t, (\mathbf{x})_t, y, z)$ can be treated along the lines of the "*Four Step Scheme*" (see, e.g., [16]). That is, one should first find a function $\mathbf{\Phi}: [0, T] \times \mathbb{C}([0, T], \mathbb{R}^n) \times \mathbb{R}^m \times \mathbb{R}^{m \times n} \mapsto \mathbb{R}^{m \times d}$ such that $\mathbf{\Phi}(t, (\mathbf{x})_t, y, z) = z\sigma(t, (\mathbf{x})_t, y, \mathbf{\Phi}(t, (\mathbf{x})_t, y, z))$, and consider $\sigma(t, (\mathbf{x})_t, y, \mathbf{\Phi}(t, (\mathbf{x})_t, y, z))$. Then we define the forward–backward martingale problem the same way as Definition 2.3 except that the functions $\hat{b}$ and $\hat{h}$ are replaced by $\hat{f}(t, (\mathbf{x})_t, y, z) = f(t, (\mathbf{x})_t, y, \mathbf{\Phi}(t, (\mathbf{x})_t, y, z))$, $f = b, h$. We leave the details to the interested reader.

We note that the Definition 2.3 looks slightly different from that of the traditional martingale problem. But one can check that they are essentially the same, modulo an application of Itô's formula. In fact, if $(\mathbb{P}, \mathbf{z})$ is a solution to the $\text{FBMP}_{x,T}(b, \sigma, h, g)$, then by Definition 2.3(i) and (iii), we have

(2.5)
$$\begin{cases} d\mathbf{x}_t = \widehat{b}(t, (\mathbf{x})_t, \mathbf{y}_t, \mathbf{z}_t)\, dt + dM_{\mathbf{x}}(t), \\ d\mathbf{y}_t = -\widehat{h}(t, (\mathbf{x})_t, \mathbf{y}_t, \mathbf{z}_t)\, dt + dM_{\mathbf{y}}(t) \\ \quad\quad = -\widehat{h}(t, (\mathbf{x})_t, \mathbf{y}_t, \mathbf{z}_t)\, dt + \mathbf{z}_t\, dM_{\mathbf{x}}(t). \end{cases}$$

Applying Itô's formula and using Definition 2.3(ii), for any $\varphi \in C^2(\mathbb{R}^n \times \mathbb{R}^m)$ and $t \in [0, T]$, one has

$$\begin{aligned} d\varphi(\mathbf{x}_t, \mathbf{y}_t) = \{&\langle \nabla_x \varphi(\mathbf{x}_t, \mathbf{y}_t), \widehat{b}(t, (\mathbf{x})_t, \mathbf{y}_t, \mathbf{z}_t)\rangle \\ &- \langle \nabla_y \varphi((\mathbf{x})_t, \mathbf{y}_t), \widehat{h}(t, (\mathbf{x})_t, \mathbf{y}_t, \mathbf{z}_t)\rangle \\ &+ \tfrac{1}{2}\operatorname{tr}\{D_{x,y}^2 \varphi(\mathbf{x}_t, \mathbf{y}_t) A(t, (\mathbf{x})_t, \mathbf{y}_t, \mathbf{z}_t)\}\} dt \\ &+ \langle \nabla_x \varphi(\mathbf{x}_t, \mathbf{y}_t), dM_{\mathbf{x}}(t)\rangle + \langle \nabla_y \varphi(\mathbf{x}_t, \mathbf{y}_t), dM_{\mathbf{y}}(t)\rangle, \end{aligned}$$

where

$$A(t, (\mathbf{x})_t, y, z) \stackrel{\triangle}{=} \begin{bmatrix} I_n \\ z \end{bmatrix} a(t, (\mathbf{x})_t, y)[I_n, z^T];$$

(2.6)

$$D^2_{x,y}\varphi = \begin{bmatrix} \partial^2_{xx}\varphi & \partial^2_{xy}\varphi \\ \partial^2_{xy}\varphi & \partial^2_{yy}\varphi \end{bmatrix}.$$

Now, if we define a differential operator

(2.7)
$$\mathcal{L}_{t,\mathbf{x},y,z} \stackrel{\triangle}{=} \tfrac{1}{2}\operatorname{tr}\{A(t,(\mathbf{x})_t,y,z)D^2_{x,y}\}$$
$$+ \langle \widehat{b}(t,(\mathbf{x})_t,y,z), \nabla_x \rangle - \langle \widehat{h}(t,(\mathbf{x})_t,y,z), \nabla_y \rangle,$$

then the fact that the $(\mathbb{P}, \mathbf{z})$ is a solution to the $\mathrm{FBMP}_{x,T}(b, \sigma, h, g)$ implies that

(2.8) $\quad C[\varphi](t) \stackrel{\triangle}{=} \varphi(\mathbf{x}_t, \mathbf{y}_t) - \varphi(x, \mathbf{y}_0) - \int_0^t \mathcal{L}_{s,(\mathbf{x})_s, \mathbf{y}_s, \mathbf{z}_s} \varphi(\mathbf{x}_s, \mathbf{y}_s)\, ds$

is a $\mathbb{P}$-martingale for all $\varphi \in C^2(\mathbb{R}^n \times \mathbb{R}^m)$. Conversely, if (2.8) is a $\mathbb{P}$-martingale for all $\varphi \in C^2(\mathbb{R}^n \times \mathbb{R}^m)$, then we can choose appropriate function $\varphi$ so that Definition 2.3 holds. In other words, Definition 2.3 actually reflects all the necessary information for a "martingale problem." But we prefer this particular form as it is more symmetric and reflects the structure of our FBSDE more explicitly.

The following theorem exhibits the connection between the weak solution and the solution to the forward–backward martingale problem.

THEOREM 2.5. *Assume $n = d$. Assume also that $\sigma = \sigma(t, (\mathbf{x})_t, y)$ is non-degenerate. Then FBSDE (2.1) has a weak solution if and only if $\mathrm{FBMP}_{x,T}(b, \sigma, h, g)$ has a solution.*

PROOF. First assume FBSDE (2.1) has a weak solution $(X, Y, Z)$ defined on a standard set-up $(\Omega, \mathcal{F}, P, \{\mathcal{F}_t\}, W)$. Note that

$$[X, Y]_t = \int_0^t \sigma(s, (X)_s, Y_s) Z_s^T \, ds.$$

Thus, since $\sigma^{-1}$ exists, we see that $Z$ is adapted to $\mathcal{F}^{X,Y}$, the filtration generated by $(X, Y)$. Using the forward equation in (2.1), we can further conclude that $W$ is also $\mathcal{F}^{X,Y}$-adapted. Therefore, without loss of generality, we may consider the canonical space $\Omega$ defined by (2.3), and let $\mathbb{P} = P \circ (X, Y)^{-1}$ be the distribution of $(X, Y)$, so that $(X, Y)$ is the canonical processes. Define $\mathbf{z}_t \stackrel{\triangle}{=} Z_t \sigma^{-1}(t, (\mathbf{x})_t, \mathbf{y}_t)$. One can check straightforwardly that $(\mathbb{P}, \mathbf{z})$ is a solution to $\mathrm{FBMP}_{x,T}(b, \sigma, h, g)$.

We next assume $\text{FBMP}_{x,T}(b,\sigma,h,g)$ has a solution $(\mathbb{P},\mathbf{z})$. Define

$$(2.9) \qquad W_t \triangleq \int_0^t \sigma^{-1}(s,(\mathbf{x})_s,\mathbf{y}_s)\, dM_{\mathbf{x}}(s).$$

Then $W$ is a continuous local martingale and $[W,W]_t = t$ by definition. Therefore, it follows from the Lévy characterization theorem (cf., e.g., [13]) we know that $W$ is a Brownian motion. Now define $Z_t \triangleq \mathbf{z}_t \sigma(t,(\mathbf{x})_t,\mathbf{y}_t)$. One can easily check that $(\mathbf{x},\mathbf{y},\mathbf{z},W)$, together with the canonical space, is a weak solution to FBSDE (2.1). $\square$

REMARK 2.6. (i) From the proof of Theorem 2.5 we see that the process $\mathbf{z}$ in Definition 2.3 is different from the martingale integrand $Z$ in FBSDE (2.1). In fact, one has the relation: $Z_t = \mathbf{z}_t \sigma(t,(\mathbf{x})_t,\mathbf{y}_t)$. Note that in the Markovian strong solution case the process $\mathbf{z}$ is actually associated directly to the gradient of the solutions to a quasilinear parabolic PDE (see, e.g., [18]).

(ii) When $\sigma$ is nondegenerate, there is an obvious one-to-one correspondence between $Z$ and $\mathbf{z}$. Thus, we shall often refer to $(\mathbb{P}, Z)$ as a solution to $\text{FBMP}_{x,T}(b,\sigma,h,g)$ as well, when the context is clear. This is particularly important in Section 5.

To conclude this section, let us give the following standing assumptions which will be used in different combinations throughout the paper:

(H1) The coefficients $(b,\sigma,h,g)$ are bounded, measurable functions, such that the mappings $(\mathbf{x},y,z) \mapsto f(t,(\mathbf{x})_t,y,z)$, $f = b, \sigma, h, g$, and $(\mathbf{x},y,z) \in \mathbb{C}([0,T];\mathbb{R}^n) \times \mathbb{R}^m \times \mathbb{R}^{m \times n}$ are uniformly continuous, uniformly in $t \in [0,T]$;
(H2) There exists a constant $K > 0$, such that $\frac{1}{K}|\lambda|^2 \leq \lambda^T \sigma\sigma^T(t,(\mathbf{x})_t,y,z)\lambda \leq K|\lambda|^2$, for all $(t,\mathbf{x},y,z) \in [0,T] \times \mathbb{C}([0,T];\mathbb{R}^n) \times \mathbb{R}^m \times \mathbb{R}^{m \times n}$ and all $\lambda \in \mathbb{R}^n$;
(H3) The mappings $t \mapsto f(t,(\mathbf{x})_t,y,z)$, $f = b, \sigma, h$, and $t \in [0,T]$ are uniformly continuous, uniformly in $(\mathbf{x},y,z) \in \mathbb{C}([0,T];\mathbb{R}^n) \times \mathbb{R}^m \times \mathbb{R}^{m \times n}$.

**3. Existence: a general result.** In this section, we study FBSDE (2.1). We note that in this section $\sigma$ may depend on $Z$. To simplify presentation in what follows, we shall assume that $\dim(X) = \dim(Y) = \dim(W) = 1$. But we note that all processes here can be higher dimensional, and all the arguments can be validated without substantial difficulties. Denoting $\|f\|_\infty = \sup|f|$ to be the usual sup-norm of a (generic) continuous function $f$, our main existence result is the following.

THEOREM 3.1. *Assume (H1), and assume that there exist a sequence of coefficients $(b_n,\sigma_n,h_n,g_n)$, $n = 1,2,\ldots$, such that:*

(i) *for $f = b, \sigma, h, g$, $\|f_n - f\|_\infty \leq \frac{1}{n}$;*

(ii) *all $(b_n, \sigma_n, h_n, g_n)$'s satisfy* (H1), *uniformly in $n$;*

(iii) *for all $n$, the FBSDE (2.1) with coefficients $(b_n, \sigma_n, h_n, g_n)$ have strong solutions $(X^n, Y^n, Z^n)$, defined on a common filtered probability space $(\Omega, \mathcal{F}, P; \mathbf{F})$ with a given $\mathbf{F}$-Brownian motion $W$;*

(iv) *denoting $Z_t^{n,\delta} \triangleq \frac{1}{\delta} \int_{(t-\delta)^+}^{t} Z_s^n \, ds$, it holds that*

$$(3.1) \qquad \lim_{\delta \to 0} \sup_n E\left\{ \int_0^T |Z_t^n - Z_t^{n,\delta}|^2 \, dt \right\} = 0.$$

*Then (2.1) admits a weak solution in the sense of Definition 2.1.*

PROOF. We shall split the proof into several steps.

*Step* 1. Denote $\Theta_t^n \triangleq ((X^n)_t, Y_t^n, Z_t^n)$ and

$$B_t^n \triangleq \int_0^t b_n(s, \Theta_s^n) \, ds; \qquad H_t^n \triangleq \int_0^t h_n(s, \Theta_s^n) \, ds; \qquad A^n(t) \triangleq \int_0^t Z_s^n \, ds;$$

$$M_t^n \triangleq \int_0^t \sigma_n(s, \Theta_s^n) \, dW_s; \qquad N_t^n \triangleq \int_0^t Z_s^n \, dW_s.$$

Consider the sequence of processes $\xi^n = (W, X^n, Y^n, B^n, H^n, A^n, M^n, N^n)$, $n = 1, 2, \ldots$, and define the canonical space $\widehat{\Omega} \triangleq \mathbb{D}([0,T])^8$ with natural filtration $\mathcal{F}$. Let $\mathbb{P}^n \triangleq P[\xi^n]^{-1} \in \mathcal{P}(\widehat{\Omega})$ be the induced probability. It is fairly easy to show that all the components in the processes $(W, X^n, Y^n, B^n, H^n, A^n, M^n, N^n)$ are quasimartingales with uniformly bounded conditional variation. For example, let $0 = t_0 < \cdots < t_m = T$ be an arbitrary partition of $[0, T]$. Then denoting $E_t \triangleq E\{\cdot | \mathcal{F}_t\}$, $t \geq 0$, one has

$$E\left\{ \sum_{i=0}^{m-1} |E_{t_i}\{Y_{t_{i+1}}^n\} - Y_{t_i}^n| + |Y_T^n| \right\}$$

$$\leq E\left\{ \sum_{i=0}^{m-1} \int_{t_i}^{t_{i+1}} |h_n(t, \Theta_t^n)| \, dt + |g_n(X_T^n)| \right\} \leq C.$$

Here and in what follows, $C > 0$ will denote a generic constant depending only on the coefficients $(b, \sigma, h, g)$ and $T$, which is allowed to vary from line to line. Thus, applying the Meyer–Zheng tightness criteria (Theorem 4 of [22]) we see that possibly along a subsequence, $\mathbb{P}^n$ converges to $\mathbb{P} \in \mathcal{P}(\widehat{\Omega})$ under the Meyer–Zheng pseudo-path topology. Consequently, $\mathbb{P}^n$ converges to $\mathbb{P}$ weakly on $\mathbb{D}([0,T])^8$, and we denote the limit to be $(W, X, Y, B, H, A, M, N)$.

*Step* 2. In the following steps, we shall identify the limit obtained in the previous step. By a slight abuse of notation, in what follows let $(W, X, Y, B, H, A, M, N)$ denote the coordinate process of $\widehat{\Omega}$. We first claim that $\mathbb{P}\{(W, X, Y,$ $B, H, A) \in \mathbb{C}([0,T])^6\} = 1$. Indeed, since by assumption (ii), all the coefficients are uniformly bounded, one can easily check that the sequence $\{(W, X^n, B^n, H^n, A^n, M^n)\}$ is tight in the space $\mathbb{C}[0,T]$ under uniform topology. [For example, if we denote $w_{M^n}(\delta) \triangleq \sup_{|s-t|\leq \delta} |M_s^n - M_t^n|$ to be the modulus of continuity of $M^n$, then it is readily seen that $E|w_{M^n}(\delta)|^2 \leq C\delta$, uniformly in $n$. Hence, by the standard tightness criteria on the space $\mathcal{P}(\mathbb{C}[0,T])$ (see, e.g., [2], Theorem 7.3), one can easily conclude that $\{M_n\}$ is tight. Other components can be argued similarly.] Consequently, the sequence $\{\mathbb{P}^n\}$ restricted to the components $(W, X, Y, B, H, A)$ converges weakly to some $\tilde{\mathbb{P}} \in \mathcal{P}(\mathbb{C}([0,T])^6)$. Since $\mathbb{C}$ is a subspace of $\mathbb{D}$, the uniqueness of the limit then leads to that $\tilde{\mathbb{P}} = \mathbb{P}|_{(W,X,Y,B,H,A)}$, proving the claim.

Next, by using the definition of weak convergence, it is fairly easy to check that

(3.2) $\quad X_t = X_0 + B_t + M_t, \qquad Y_t = Y_0 - H_t + N_t \qquad \forall t \in [0,T), \ \mathbb{P}\text{-a.s.}$

Clearly, under probability $\mathbb{P}$, $W$ is a Brownian motion. Since $X, B, H, M$ are all continuous, noting that $\sup_n E \int_0^T |Z_t^n|^2 \, dt < \infty$, it follows from [22], Theorem 11, that $M$, $N$ are both martingales. Further, applying [22], Theorem 10, we conclude that $A$ is absolutely continuous, $\mathbb{P}$-a.s.; and $A_t = \int_0^t Z_s \, ds$ with $E^{\mathbb{P}} \int_0^T |Z_t|^2 \, dt < \infty$.

*Step* 3. We show that $B_t = \int_0^t b(s, \Theta_s) \, ds$ and $H_t = \int_0^t h(s, \Theta_s) \, ds$, $\forall t$, $\mathbb{P}$-a.s. To this end, we note that the function $b$ is uniformly continuous on $z$. Thus, for any $\varepsilon > 0$, there exists $\varepsilon_0 > 0$ so that $|b(t, (\mathbf{x})_t, y, z_1) - b(t, (\mathbf{x})_t, y, z_2)| \leq \varepsilon$ whenever $|z_1 - z_2| \leq \varepsilon_0$. Furthermore, (3.1), we can choose $\delta_0 > 0$ such that for any $\delta \leq \delta_0$ it holds that

(3.3) $$\sup_n E\left\{\int_0^T |Z_t^n - Z_t^{n,\delta}|^2 \, dt\right\} \leq \varepsilon \varepsilon_0^2.$$

Now let us denote $Z_t^\delta \triangleq \frac{1}{\delta}[A_t - A_{t-\delta}]$, where $A_t \triangleq 0$ for $t < 0$. Then by assumption (i) and the definition of $\{\mathbb{P}^n\}$, one verifies easily that

$$E^{\mathbb{P}}\left\{\left|B_t - \int_0^t b(s, \Theta_s) \, ds\right|\right\}$$

$$= \lim_{\delta \to 0} E^{\mathbb{P}}\left\{\left|B_t - \int_0^t b(s, (X)_s, Y_s, Z_s^\delta) \, ds\right|\right\}$$

(3.4) $\quad = \lim_{\delta \to 0} \lim_{n \to \infty} E^{\mathbb{P}^n}\left\{\left|B_t - \int_0^t b(s, (X)_s, Y_s, Z_s^\delta) \, ds\right|\right\}$

$$= \lim_{\delta \to 0} \lim_{n \to \infty} E\left\{\left|\int_0^t b_n(s, \Theta_s^n) \, ds - \int_0^t b(s, (X^n)_s, Y_s^n, Z_s^{n,\delta}) \, ds\right|\right\}$$

$$\leq \lim_{\delta \to 0} \lim_{n \to \infty} E\left\{\int_0^T |b(s, (X^n)_s, Y_s^n, Z_s^n) - b(s, (X^n)_s, Y_s^n, Z_s^{n,\delta})| \, ds\right\}.$$

Furthermore, denote $\Delta b_s^{n,\delta} \triangleq b(s, (X^n)_s, Y_s^n, Z_s^n) - b(s, (X^n)_s, Y_s^n, Z_s^{n,\delta})$, for $s \in [0, T]$, using (3.3) and the boundedness of $b$ we deduce that

$$E \int_0^T |\Delta b_s^{n,\delta}| \, ds$$

$$= E \int_0^T \{|\Delta b_s^{n,\delta}|[1_{\{|Z_s^n - Z_s^{n,\delta}| \leq \varepsilon_0\}} + 1_{\{|Z_s^n - Z_s^{n,\delta}| > \varepsilon_0\}}] \, ds\}$$

(3.5)
$$\leq T\varepsilon + CE\left\{\int_0^T 1_{\{|Z_s^n - Z_s^{n,\delta}| > \varepsilon_0\}} \, ds\right\}$$

$$\leq T\varepsilon + \frac{C}{\varepsilon_0^2} E\left\{\int_0^T |Z_s^n - Z_s^{n,\delta}|^2 \, ds\right\}$$

$$\leq (T + C)\varepsilon.$$

Since $\varepsilon$ is arbitrary, we get $E^{\mathbb{P}}\{|B_t - \int_0^t b(s, \Theta_s) \, ds|\} = 0$. An almost identical proof also shows that $E^{\mathbb{P}}\{|H_t - \int_0^t h(s, \Theta_s) \, ds|\} = 0$. This completes the claim.

*Step* 4. We now show that $N_t = \int_0^t Z_s \, dW_s$, $\forall t < T$, $\mathbb{P}$-a.s. To see this, first note that $N$ is càdlàg and $\int_0^t Z_s \, dW_s$ is continuous, then it suffices to show that

(3.6)
$$I \triangleq E^{\mathbb{P}}\left\{\int_0^T \left|N_t - \int_0^t Z_s \, dW_s\right|^2 dt\right\} = 0.$$

Again, we shall use the fact that $\mathbb{P}^n \to \mathbb{P}$ weakly. But in this case, we should note that in (3.6) the stochastic integral is generally unbounded, therefore, an extra truncation procedure is necessary. Indeed, applying the Monotone Convergence Theorem, we see that to prove (3.6) it suffices to show that

(3.7) $$I_R \triangleq E^{\mathbb{P}}\left\{\left[\int_0^T \left|N_t - \int_0^t Z_s \, dW_s\right|^2 dt\right] \wedge R\right\} = 0 \qquad \forall R > 0.$$

We now fix $R > 0$, and notice the following simple fact:

$$(a + b) \wedge R \leq a \wedge R + b \leq a + b \qquad \forall a, b \geq 0.$$

By definition of $Z^\delta$ one checks that

(3.8) $$\lim_{\delta \to 0} E^{\mathbb{P}}\left\{\int_0^T |Z_t - Z_t^\delta|^2 \, dt\right\} = 0.$$

But this implies that

$$I_R \leq \varlimsup_{\delta \to 0} E^{\mathbb{P}}\left\{\left[\int_0^T \left|N_t - \int_0^t Z_s^\delta \, dW_s\right|^2 dt\right] \wedge R\right\} \triangleq \varlimsup_{\delta \to 0} I^\delta.$$

Let $\pi: 0 = t_0 < \cdots < t_m = T$ be any partition of $[0,T]$. Then

$$I^\delta \leq CE^{\mathbb{P}}\left\{\left[\sum_{j=0}^{m-1}\int_{t_j}^{t_{j+1}}\left|N_t - \sum_{i=0}^{j-1} Z_{t_i}^\delta[W_{t_{i+1}} - W_{t_i}]\right|^2 dt\right] \wedge R\right\}$$

(3.9)
$$+ \frac{C}{\delta^2} E^{\mathbb{P}}\{I^{\pi,\delta}\},$$

where

$$I^{\pi,\delta} \triangleq \sum_{j=0}^{m-1}\int_{t_j}^{t_{j+1}}\left|\sum_{i=0}^{j-1}\int_{t_i}^{t_{i+1}}[(A_s - A_{s-\delta}) - (A_{t_i} - A_{t_i-\delta})]\,dW_s\right.$$

$$\left.+ \int_{t_j}^{t}[A_s - A_{s-\delta}]\,dW_s\right|^2 dt.$$

A similar calculation shows further that (changing $\mathbb{P}$ to $\mathbb{P}^n$ when necessary),

$$E^{\mathbb{P}}\left\{\left[\sum_{j=0}^{m-1}\int_{t_j}^{t_{j+1}}\left|N_t - \sum_{i=0}^{j-1} Z_{t_i}^\delta[W_{t_{i+1}} - W_{t_i}]\right|^2 dt\right] \wedge R\right\}$$

$$= \lim_{n \to \infty} E\left\{\left[\sum_{j=0}^{m-1}\int_{t_j}^{t_{j+1}}\left|\int_0^t Z_s^n\,dW_s\right.\right.\right.$$

(3.10)
$$\left.\left.\left.- \sum_{i=0}^{j-1} Z_{t_i}^{n,\delta}[W_{t_{i+1}} - W_{t_i}]\right|^2 dt\right] \wedge R\right\}$$

$$\leq C\,\overline{\lim}_{n \to \infty} E\left\{\int_0^T \left|\int_0^t Z_s^n\,dW_s - \int_0^t Z_s^{n,\delta}\,dW_s\right|^2 dt\right\}$$

$$+ C\,\overline{\lim}_{n \to \infty} \frac{E^{\mathbb{P}^n}\{I^{\pi,\delta}\}}{\delta^2}.$$

Now, let us denote $\Delta A_t^s = A_t - A_s = \int_s^t Z_r\,dr$, $0 \leq s \leq t \leq T$. We see that

$$E^{\mathbb{P}}\{I^{\pi,\delta}\} = E^{\mathbb{P}}\left\{\sum_{j=0}^{m-1}\int_{t_j}^{t_{j+1}}\left[\sum_{i=0}^{j-1}\int_{t_i}^{t_{i+1}} |\Delta A_s^{t_i} - \Delta A_{s-\delta}^{t_i-\delta}|^2\,ds\right.\right.$$

$$\left.\left.+ \int_{t_j}^{t}|\Delta A_s^{s-\delta}|^2\,ds\right]dt\right\}$$

(3.11)
$$\leq CE^{\mathbb{P}}\left\{\sum_{j=0}^{m-1}\int_{t_j}^{t_{j+1}}\left[\sum_{i=0}^{j-1}\int_{t_i}^{t_{i+1}} (s-t_i)\int_0^T |Z_r|^2\,dr\,ds\right.\right.$$

$$\left.\left.+ \delta \int_{t_j}^{t}\int_{s-\delta}^{s}|Z_r|^2\,dr\,ds\right]dt\right\}$$

$$\le C|\pi| E^{\mathbb{P}}\Big\{\int_0^T |Z_t|^2\,dt\Big\} \le C|\pi|.$$

Similarly, one shows that $E^{\mathbb{P}^n}\{I^{\pi,\delta}\} \le C|\pi|$. This, together with (3.10) and (3.11), reduces (3.9) to the following:

$$I^\delta \le C\varlimsup_n E\Big\{\int_0^T \Big|\int_0^t Z_s^n\,dW_s - \int_0^t Z_s^{n,\delta}\,dW_s\Big|^2\,dt\Big\} + \frac{C|\pi|}{\delta^2}.$$

Since $\pi$ is arbitrary, we have

$$I^\delta \le C\varlimsup_n E\Big\{\int_0^T \Big|\int_0^t Z_s^n\,dW_s - \int_0^t Z_s^{n,\delta}\,dW_s\Big|^2\,dt\Big\}$$

$$= C\varlimsup_n E\Big\{\int_0^T \int_0^t |Z_s^n - Z_s^{n,\delta}|^2\,dt\Big\}$$

$$\le C\varlimsup_n E\Big\{\int_0^T |Z_t^n - Z_t^{n,\delta}|^2\,dt\Big\}.$$

Now, applying (3.1) we prove (3.7).

*Step* 5. We next show that $Y_T = g(X_T)$. We should note that in the last step we actually proved that the process $N$, whence $Y$, is *continuous* on $[0,T)$. Therefore, by defining $N_T = N_{T-}$, we can assume that $Y$ is (left) continuous at $T$ as well. Thus, in what follows, we shall only check that $Y_{T-} = \lim_{\varepsilon\downarrow 0} \frac{1}{\varepsilon}\int_{T-\varepsilon}^T Y_s\,ds = g(X_T)$, $\mathbb{P}$-a.s. To this end, we note that for any $\varepsilon > 0$, one has

$$E^{\mathbb{P}}\Big\{\Big|\frac{1}{\varepsilon}\int_{T-\varepsilon}^T Y_t\,dt - g(X_T)\Big|^2\Big\}$$

$$= \lim_{n\to\infty} E^{\mathbb{P}^n}\Big\{\Big|\frac{1}{\varepsilon}\int_{T-\varepsilon}^T Y_t\,dt - g(X_T)\Big|^2\Big\}$$

$$= \lim_{n\to\infty} E\Big\{\Big|\frac{1}{\varepsilon}\int_{T-\varepsilon}^T Y_t^n\,dt - Y_T^n\Big|^2\Big\}$$

$$= \lim_{n\to\infty} E\Big\{\Big|\frac{1}{\varepsilon}\int_{T-\varepsilon}^T \Big[\int_t^T h_n(s,\Theta_s^n)\,ds - \int_t^T Z_s^n\,dW_s\Big]\,dt\Big|^2\Big\}$$

$$\le \varlimsup_{n\to\infty} E\Big\{\int_{T-\varepsilon}^T |Z_t^n|^2\,dt\Big\} + C\varepsilon$$

$$\le \varlimsup_{n\to\infty} 2E\Big\{\int_{T-\varepsilon}^T [|Z_t^n - Z_t^{n,\delta}|^2 + |Z_t^{n,\delta}|^2]\,dt\Big\} + C\varepsilon.$$

We should point out that unlike step 4, in the above we do not need to apply the truncation technique, thanks to the boundedness of both process $Y$ and function $g$. Furthermore, following the arguments of step 4, we fix $\delta > 0$ and let $\pi : T - \varepsilon = t_0 < \cdots < t_m = T$ be any partition of $[T - \varepsilon, T]$. Again, we denote for any process $\xi$, $\Delta \xi_t^s = \xi_t - \xi_s$, $0 \le s \le t \le T$. Then by definition of $Z^{n,\delta}$, we have

$$
\begin{aligned}
&E\left\{\int_{T-\varepsilon}^{T} |Z_t^{n,\delta}|^2 \, dt\right\} \\
&= E\left\{\left|\int_{T-\varepsilon}^{T} \frac{1}{\delta}[A_t^n - A_{t-\delta}^n] \, dW_t\right|^2\right\} \\
&= E\left\{\left|\sum_{j=0}^{m-1} \int_{t_j}^{t_{j+1}} [\Delta[A^n]_{t_j}^{t_j-\delta} + \Delta[A^n]_t^{t_j} - \Delta[A^n]_{t-\delta}^{t_j-\delta}] \, dW_t\right|^2\right\}
\end{aligned}
$$
(3.12)
$$
\begin{aligned}
&\le CE\left\{\left|\sum_{j=0}^{m-1} \Delta[A^n]_{t_j}^{t_j-\delta} \Delta W_{t_{j+1}}^{t_j}\right|^2 \right. \\
&\qquad \left. + \sum_{j=0}^{m-1} \int_{t_j}^{t_{j+1}} [|\Delta[A^n]_t^{t_j}|^2 + |\Delta[A^n]_{t-\delta}^{t_j-\delta}|^2] \, dt\right\} \\
&\le CE^{\mathbb{P}^n}\left\{\left|\sum_{j=0}^{m-1} \Delta A_{t_j}^{t_j-\delta} \Delta W_{t_{j+1}}^{t_j}\right|^2\right\} + \frac{C|\pi|}{\delta^2},
\end{aligned}
$$

where the last inequality is due to a similar argument for (3.11). By the weak convergence of $\mathbb{P}^n$ and by using the above arguments in a reverse order, we see that there exists $N$ such that, for any $n > N$,

(3.13)
$$
\begin{aligned}
E^{\mathbb{P}^n}&\left\{\left|\sum_{j=0}^{m-1} \Delta A_{t_j}^{t_j-\delta} \Delta W_{t_{j+1}}^{t_j}\right|^2\right\} \\
&\le E^{\mathbb{P}}\left\{\left|\sum_{j=0}^{m-1} \Delta A_{t_j}^{t_j-\delta} \Delta W_{t_{j+1}}^{t_j}\right|^2\right\} + \varepsilon \\
&\le E^{\mathbb{P}}\left\{\int_{T-\varepsilon}^{T} |Z_t^{\delta}|^2 \, dt\right\} + \frac{C|\pi|}{\delta^2} + \varepsilon \\
&\le 2E^{\mathbb{P}}\left\{\int_{T-\varepsilon}^{T} [|Z_t|^2 + |Z_t - Z_t^{\delta}|^2] \, dt\right\} + \frac{C|\pi|}{\delta^2} + \varepsilon.
\end{aligned}
$$

Combining (3.12) and (3.13), we obtain that

$$
\begin{aligned}
E^{\mathbb{P}}&\left\{\left|\frac{1}{\varepsilon}\int_{T-\varepsilon}^{T} Y_t \, dt - g(X_T)\right|^2\right\} \\
&\le C\left[\overline{\lim_{n}} E\left\{\int_{T-\varepsilon}^{T} |Z_t^n - Z_t^{n,\delta}|^2 \, dt\right\}\right.
\end{aligned}
$$

$$+ E^{\mathbb{P}}\bigg\{\int_{T-\varepsilon}^{T}[|Z_t|^2 + |Z_t - Z_t^\delta|^2]\,dt\bigg\} + \frac{|\pi|}{\delta^2} + \varepsilon\bigg].$$

First letting $|\pi| \to 0$, then letting $\delta \to 0$, finally letting $\varepsilon \to 0$, and applying Fatou's Lemma we derive $E^{\mathbb{P}}\{|Y_{T-} - g(X_T)|^2\} = 0$.

*Step* 6. Finally, we note that $X^n$, $Y^n$'s all have better regularity than $Z^n$'s. Following the same arguments in step 4, we can show that

$$M_t = \int_0^t \sigma(s,\Theta_s)\,dW_s \qquad \forall t \in [0,T],\ \mathbb{P}\text{-a.s.}$$

The proof is now complete. $\square$

**4. The Markovian case.** In this section, we further explore our main existence result Theorem 3.1. It is clear that the key condition in that theorem is the assumption (3.1) on the martingale integrands $\{Z^n\}$, which in a sense represents the "path regularity" of the sequence $\{Z^n\}$ or as a certain "tightness" criterion in the space $L^2$. Without digging deep on this issue, in this section, we shall investigate some special cases where condition (3.1) is satisfied. To begin with, let us consider the following *Markovian FBSDE*:

$$(4.1) \qquad \begin{cases} X_t = x + \int_0^t \sigma(s, X_s, Y_s)\,dW_s; \\ Y_t = g(X_T) + \int_t^T h(s, X_s, Y_s, Z_s)\,ds - \int_t^T Z_s\,dW_s; \end{cases}$$

where all processes $X$, $Y$, $Z$ and $W$ are one-dimensional. We note that the assumption that the drift of the forward equation $b \equiv 0$, is merely for simplicity. In fact, the case when $b \neq 0$ can be easily deduced to such a form via a standard Girsanov transformation, especially in the case when $\sigma$ is nondegenerate and $h$ is allowed to have *linear growth* in $Z$. The assumption that all processes (especially $X$) is one dimensional is more technical, since we are going to apply a result by Nash [23] in Lemma 4.2 below. We believe that these restrictions can all be removed with more technicalities, and we shall leave them to our future publications, as these are not the main points of this paper.

The Markovian nature of the FBSDE now enables us to apply the idea of the Four Step Scheme initiated in [16]. To be more precise, we shall look for solutions to (4.1) for which the relation $Y_t = u(t, X_t)$ holds, where $u$ is a viscosity solution to the PDE

$$(4.2) \qquad \begin{cases} u_t + \frac{1}{2}\sigma^2(t,x,u)u_{xx} + h(t,x,u,u_x\sigma) = 0; \\ u(T,x) = g(x). \end{cases}$$

In [21], we proved that if besides the standing assumptions (H1) and (H2), the coefficients of (4.1) are uniformly Hölder continuous and the comparison theorem for the viscosity solution to PDE (4.2) holds true, then the weak solution to FBSDE (4.1) exists and is unique in law, and that $Y_t = u(t, X_t)$ where $u$ is the unique viscosity to PDE (4.2). The proof there relied heavily on some a priori gradient estimates of the solution $u$ of the PDE (4.2). However, these estimates are no longer valid under merely the assumptions (H1) and (H2), we shall therefore turn to Theorem 3.1.

We first recall a result which is standard in the literature (see, e.g., [7]).

LEMMA 4.1. *Assume* (H1) *and* (H2). *There exist a viscosity solution $u$ to PDE (4.2) and constants $C$ and $\alpha > 0$ such that for any $x, y \in \mathbb{R}$ and any $0 \le s \le t < T$,*

$$(4.3) \qquad |u(s,x) - u(t,y)| \le \frac{C}{(T-t)^{\alpha/2}}[|x-y|^\alpha + |t-s|^{\alpha/2}].$$

We next establish an a priori estimate for the following linear PDE:

$$(4.4) \qquad \begin{cases} u_t + \frac{1}{2}\sigma^2(t,x)u_{xx} = 0; \\ u(T,x) = g(x). \end{cases}$$

LEMMA 4.2. *Assume that $\sigma$ is smooth in $x$ and* (H2) *holds, and that $g \in C^2$ with $\|g\|_\infty + \|g'\|_\infty + \|g''\|_\infty \le K$. Then there exist constants $C$ and $\alpha > 0$, depending only on $K$ and $T$, and are independent of the derivatives of $\sigma$, such that*

$$(4.5) \quad \|u_x\|_\infty \le C; \qquad |u_x(s,x) - u_x(t,y)| \le C[|s-t|^{\alpha/2} + |x-y|^\alpha].$$

PROOF. Note that $v \stackrel{\triangle}{=} u_x$ is the solution to the following PDE in divergence form:

$$\begin{cases} v_t + \frac{1}{2}(\sigma^2 v_x)_x = 0; \\ v(T,x) = g'(x). \end{cases}$$

The result follows from some well-known results of Nash [23]. □

Let us now consider the FBSDE corresponding to the simplified PDE (4.4):

$$(4.6) \qquad \begin{cases} X_t = x + \int_0^t \sigma(s, X_s)\, dW_s; \\ Y_t = g(X_T) - \int_t^T Z_s\, dW_s; \end{cases} \quad t \in [0, T].$$

LEMMA 4.3. *Assume that $\sigma$ satisfies* (H1) *and* (H2), *and assume that there exists a sequence of functions $\{\sigma_n(t,x)\}$ such that:*

(i) *each $\sigma_n$ is smooth in $x$ and satisfy (H2), such that $\|\sigma'_n\|_\infty \le C_n$, for all $n$;*

(ii) *$\sigma_n \to \sigma$ uniformly;*

(iii) *$g_n \in C^2$, with $\|g_n\|_\infty + \|g'_n\|_\infty + \|g''_n\|_\infty \le C$, for some generic constant $C > 0$.*

For each $n$, let $(X^n, Y^n, Z^n)$ be the strong solution to the FBSDE (4.6) with coefficients $\sigma_n$ and $g_n$. Then denoting $Z^n_t = 0$ for $t < 0$, it holds that

$$(4.7) \qquad \lim_{\delta \to 0} \sup_n E\left\{\int_0^T |Z^n_t - Z^n_{t-\delta}|^2 \, dt\right\} = 0.$$

PROOF. First, by the Four Step Scheme we have $Z^n_t = u^n_x \sigma_n(t, X^n_t)$, where $u^n$ is the classical solution to the following PDE:

$$(4.8) \qquad \begin{cases} u^n_t + \frac{1}{2}\sigma_n^2(t,x) u^n_{xx} = 0; \\ u^n(T,x) = g_n(x). \end{cases}$$

We shall assume from now on that $\xi_t = 0$, $t < 0$ for all processes $\xi \in L^2([0,T] \times \Omega)$. Applying Lemma 4.2, we see that $u^n$'s satisfy (4.5) uniformly (in $n$). Thus, denoting $C > 0$ to be all the generic constant, we have

$$E\left\{\int_0^T |Z^n_t - Z^n_{t-\delta}|^2 \, dt\right\}$$
$$\le CE\left\{\int_0^T [|u^n_x(t, X^n_t) - u^n_x(t-\delta, X^n_{t-\delta})|^2 \right.$$
$$\left. + |\sigma_n(t, X^n_t) - \sigma_n(t-\delta, X^n_{t-\delta})|^2] \, dt\right\}$$
$$\le CE\left\{\int_0^T [\delta^\alpha + |X^n_t - X^n_{t-\delta}|^{2\alpha} + |\sigma_n(t, X^n_t) - \sigma_n(t-\delta, X^n_{t-\delta})|^2] \, dt\right\}$$
$$\le C\delta^\alpha + CE\left\{\int_0^T |\sigma_n(t, X^n_t) - \sigma_n(t-\delta, X^n_{t-\delta})|^2 \, dt\right\}.$$

Therefore, it suffices to show that

$$(4.9) \qquad \lim_{\delta \to 0} \sup_n E\left\{\int_0^T |\sigma_n(t, X^n_t) - \sigma_n(t-\delta, X^n_{t-\delta})|^2 \, dt\right\} = 0.$$

We should note that if $\sigma_n$ is uniformly continuous in $t$, then (4.9) is obviously true. But under (H1), $\sigma$ may not even be continuous (!). Therefore, we shall use (ii) instead.

To this end, first note that by approximation using processes with continuous paths if necessary, one can show that, for any process $\xi \in L^2([0,T] \times \Omega)$,

$$(4.10) \qquad \lim_{\delta \to 0} E\left\{\int_0^T |\xi_t - \xi_{t-\delta}|^2 \, dt\right\} = 0.$$

Next, by assumption (ii), we see that for any $\varepsilon > 0$, there exists $N_0 > 0$, such that
$$|\sigma_n(t,x) - \sigma(t,x)| < \varepsilon \qquad \forall (t,x) \in [0,T] \times \mathbb{R}^d,$$
whenever $n > N_0$. Thus, for $n > N_0$, we have

(4.11)
$$E\left\{\int_0^T |\sigma_n(t, X_t^n) - \sigma_n(t-\delta, X_{t-\delta}^n)|^2 \, dt\right\}$$
$$\leq CE\left\{\int_0^T |\sigma(t, X_t^n) - \sigma(t-\delta, X_{t-\delta}^n)|^2 \, dt\right\} + C\varepsilon^2.$$

Furthermore, note that the distributions of the sequence $\{X^n\}$ are obviously tight (see the previous section), and it is readily seen that $X^n$ must converge to $X$ in distribution, where $X$ is the unique weak solution to the following SDE:
$$X_t = x + \int_0^t \sigma(s, X_s) \, dW_s.$$

By the Skorohod representation theorem, we can assume without loss of generality that on a common probability space, still denote it by $(\Omega, \mathcal{F}, P)$, $X^n$ converges to $X$ almost surely. Now applying the Bounded Convergence Theorem, and changing $N_0$ if necessary, we can modify (4.11) to the following:

(4.12)
$$E\left\{\int_0^T |\sigma_n(t, X_t^n) - \sigma_n(t-\delta, X_{t-\delta}^n)|^2 \, dt\right\}$$
$$\leq CE\left\{\int_0^T |\sigma(t, X_t) - \sigma(t-\delta, X_{t-\delta})|^2 \, dt\right\} + C\varepsilon^2 \qquad \forall n > N_0.$$

Finally, denote
$$\begin{cases} I_n(\delta) \triangleq E\left\{\int_0^T |\sigma_n(t, X_t^n) - \sigma_n(t-\delta, X_{t-\delta}^n)| \, dt\right\}, & n = 1, 2, \ldots, \\ I(\delta) \triangleq E\left\{\int_0^T |\sigma(t, X_t) - \sigma(t-\delta, X_{t-\delta})| \, dt\right\}, \end{cases}$$

(4.12) then leads to
$$\sup_n I_n(\delta) \leq \sup_{n \leq N_0} I_n(\delta) + \sup_{n > N_0} I_n(\delta) \leq \sum_{n=1}^{N_0} I_n(\delta) + CI(\delta) + C\varepsilon^2.$$

Consequently, we see that (4.9) follows by first letting $\delta \to 0$ and applying (4.10) in the above, and then letting $\varepsilon \to 0$. The proof is complete. □

We note that in Lemma 4.3 the assumptions on $g_n$ are rather strong. The following lemma is a weaker alternative.

LEMMA 4.4. *Assume* (H1), (H2) *as well as* (i) *and* (ii) *in Lemma* 4.3. *Assume further that* $\|g_n\|_\infty \leq K$, *and for each* $n$, *there exists constant* $C_n > 0$ *such that* $\|g'_n\|_\infty$, $\|g''_n\|_\infty \leq C_n$. *Again, denote* $(X^n, Y^n, Z^n)$ *be the strong solutions to FBSDE* (4.6) *with coefficients* $\sigma_n$ *and* $g_n$. *Then the following conclusions hold:*

(i) *If* $g_n$'s *are uniformly continuous, uniformly on* $n$, *then* (4.7) *holds true.*

(ii) *In general, for any* $\varepsilon > 0$,

(4.13) $$\lim_{\delta \to 0} \sup_n E\left\{ \int_0^{T-\varepsilon} |Z_t^n - Z_{t-\delta}^n|^2 \, dt \right\} = 0.$$

PROOF. (i) Since $g_n$ is uniformly continuous, for any $\varepsilon > 0$, we may find a mollifier of $g_n$, denoted by $\bar{g}_n$, such that $\|\bar{g}_n - g_n\|_\infty \leq \sqrt{\varepsilon}$; and that

$$\|\bar{g}_n\|_\infty \leq K, \qquad \|\bar{g}'_n\|_\infty \leq C_\varepsilon, \qquad \|\bar{g}''_n\|_\infty \leq C_\varepsilon.$$

We should note that since the uniform continuity of $g_n$'s is assumed to be *uniform* in $n$, it follows from the standard mollification procedure that the constant $C_\varepsilon$ can be chosen to be independent of $n$ as well. In particular, we have $\|\bar{g}_n(X_T^n) - g_n(X_T^n)\|_{L^2(\Omega)}^2 \leq \varepsilon$. Now, let $(\bar{Y}^n, \bar{Z}^n)$ be the solution to BSDE

$$\bar{Y}_t^n = \bar{g}_n(X_T^n) - \int_t^T \bar{Z}_s^n \, dW_s.$$

Then by the standard estimates for BSDEs we see that $E\{\int_0^T |\bar{Z}_t^n - Z_t^n|^2 \, dt\} \leq \varepsilon$. Applying Lemma 4.3, we have

$$\lim_{\delta \to 0} \sup_n E\left\{ \int_0^T |\bar{Z}_t^n - \bar{Z}_{t-\delta}^n|^2 \, dt \right\} = 0.$$

Thus,

$$\varlimsup_{\delta \to 0} \sup_n E\left\{ \int_0^T |Z_t^n - Z_{t-\delta}^n|^2 \, dt \right\}$$
$$\leq C \varlimsup_{\delta \to 0} \sup_n E\left\{ \int_0^T [|Z_t^n - \bar{Z}_t^n|^2 + |\bar{Z}_t^n - \bar{Z}_{t-\delta}^n|^2 + |\bar{Z}_{t-\delta}^n - Z_{t-\delta}^n|^2] \, dt \right\}$$
$$\leq C\varepsilon.$$

Since $\varepsilon$ is arbitrary, the result follows.

(ii) In this case, we let $u_n$ be the classical solution to the PDE (4.8). By Lemma 4.1, we know that $u_n(T - \varepsilon, \cdot)$ is uniformly Hölder-$\alpha$ continuous in $x$. Note that the processes $(X^n, Y^n, Z^n)$ also satisfy

$$\begin{cases} X_t^n = x + \int_0^t \sigma_n(s, X_s^n) \, dW_s; \\ Y_t^n = u_n(T - \varepsilon, X_{T-\varepsilon}^n) - \int_t^{T-\varepsilon} Z_s^n \, dW_s. \end{cases}$$

Applying part (i) to the equation above, we obtain the result. □

We are now ready to prove the main result of this section.

THEOREM 4.5. *Assume* (H1) *and* (H2). *Then FBSDE (4.1) admits a weak solution* $(\Omega, \mathcal{F}, \mathbb{P}; \mathbf{F}, W, X, Y, Z)$. *Moreover, it holds that* $Y_t = u(t, X_t)$, $t \in [0, T]$, $\mathbb{P}$-*a.s., where* $u$ *is a viscosity solution to PDE (4.2) satisfying (4.3).*

PROOF. Let $\sigma_n, h_n, g_n$ be the standard molifiers of $\sigma, h, g$, such that $g_n$ is uniformly continuous and

$$\|\sigma_n\|_\infty + \|h_n\|_\infty + \|g_n\|_\infty \le K; \quad \frac{1}{K} \le \sigma_n \le K.$$

Let $(X^n, Y^n, Z^n)$ be the strong solution to the following FBSDE:

$$(4.14) \quad \begin{cases} X_t^n = x + \int_0^t \sigma_n(s, X_s^n, Y_s^n) \, dW_s; \\ Y_t^n = g_n(X_T^n) + \int_t^T h_n(s, X_s^n, Y_s^n, Z_s^n) \, ds - \int_t^T Z_s^n \, dW_s, \end{cases}$$

and $u^n$ the classical solution to the following PDE:

$$(4.15) \quad \begin{cases} u_t^n + \frac{1}{2} \sigma_n^2(t, x, u^n(t, x)) u_{xx}^n \\ \quad + h_n(t, x, u^n(t, x), u_x^n \sigma(t, x, u^n(t, x))) = 0; \\ u^n(T, x) = g_n(x). \end{cases}$$

Denote $\tilde{\sigma}_n(t, x) \triangleq \sigma_n(t, x, u^n(t, x))$ and $\tilde{h}_n(t, x) \triangleq h_n(t, x, u^n(t, x), u_x^n \sigma(t, x, u^n(t, x)))$. Note that the solutions $u^n$'s are bounded and uniformly continuous in $(t, x)$ (actually uniformly Hölder continuous in $(t, x)$ if $g_n$'s are uniformly Hölder continuous, cf. [7]). Thus, applying Arzelà–Ascoli theorem, it follows that $u^n \to u$ uniformly on compact sets, where $u$ is the unique viscosity solution to PDE (4.2) and $u$ is also uniformly continuous in $(t, x)$. Therefore, we have

$$(4.16) \quad \tilde{\sigma}_n(t, x) \to \tilde{\sigma}(t, x) \triangleq \sigma(t, x, u(t, x)),$$

and the limit is uniform as well. Note that

$$Y_t^n = u_n(t, X_t^n), \quad Z_t^n = u_x^n(t, X_t^n) \sigma_n(t, X_t^n, u_n(t, X_t^n)).$$

We have

$$\begin{cases} X_t^n = x + \int_0^t \tilde{\sigma}_n(s, X_s^n) \, dW_s; \\ Y_t^n = g_n(X_T^n) + \int_t^T \tilde{h}_n(s, X_s^n) \, ds - \int_t^T Z_s^n \, dW_s. \end{cases}$$

Let us decompose the FBSDE above into the following two BSDEs: for $0 \leq t \leq s \leq T$,

(4.17)
$$Y_t^{n,\infty} = g_n(X_T^n) - \int_t^T Z_r^{n,\infty} \, dW_r;$$

$$Y_t^{n,s} = \tilde{h}_n(s, X_s^n) - \int_t^s Z_r^{n,s} \, dW_r.$$

We should point out here that the family of processes $\{Z_r^{n,s} : 0 \leq r \leq s \leq T\}$ is actually a random field defined on $[0,T]^2$, restricted to the triangular domain $0 \leq r \leq s \leq T$, such that for each $s \in [0,T]$ and $r \in [0,s]$, $Z_r^{n,s}$ is $\mathcal{F}_r$-measurable. (In fact, $Z^{n,s}$ has the representation: $Z_r^{n,s} = E\{\tilde{h}_n(s, X_s^n) \nabla X_s^n | \mathcal{F}_r\} [\nabla X_r^n]^{-1}$, cf. [19].) Furthermore, a simple computation using Fubini's theorem shows that

$$Y_t^n = Y_t^{n,\infty} + \int_t^T Y_t^{n,s} \, ds; \qquad Z_t^n = Z_t^{n,\infty} + \int_t^T Z_t^{n,s} \, ds, \qquad t \in [0,T].$$

Now, for any $\varepsilon > 0$ and $s \in [0,T]$, we apply Lemma 4.4(i) and (ii), respectively, to get

(4.18)
$$\limsup_{\delta \to 0} E\left\{ \int_0^T |Z_t^{n,\infty} - Z_{t-\delta}^{n,\infty}|^2 \, dt + \int_0^{s-\varepsilon} |Z_t^{n,s} - Z_{t-\delta}^{n,s}|^2 \, dt \right\} = 0.$$

Moreover, denoting $C > 0$ to be a generic constant, depending only on the bounds of the coefficients and $T$, and allowed to vary from line to line, we have

$$E\left\{ \int_0^s |Z_t^{n,s} - Z_{t-\delta}^{n,s}|^2 \, dt \right\} \leq CE\left\{ \int_0^s |Z_t^{n,s}|^2 \, dt \right\} \leq C,$$

thanks to (H1). Hence, for any $t \geq 0$, $\varepsilon > 0$, with $t + \varepsilon \leq T$, one has

$$E\left\{ \int_0^T |Z_t^n - Z_{t-\delta}^n|^2 \, dt \right\}$$

$$\leq CE\left\{ \int_0^T \left[ |Z_t^{n,\infty} - Z_{t-\delta}^{n,\infty}| + \int_{t+\varepsilon}^T |Z_t^{n,s} - Z_{t-\delta}^{n,s}| \, ds \right.\right.$$

$$\left.\left. + \int_t^{t+\varepsilon} |Z_t^{n,s} - Z_{t-\delta}^{n,s}| \, ds + \int_{t-\delta}^t |Z_{t-\delta}^{n,s}| \, ds \right]^2 dt \right\}$$

$$\leq CE\left\{ \int_0^T \left[ |Z_t^{n,\infty} - Z_{t-\delta}^{n,\infty}|^2 + \int_{t+\varepsilon}^T |Z_t^{n,s} - Z_{t-\delta}^{n,s}|^2 \, ds \right.\right.$$

$$\left.\left. + \varepsilon \int_t^{t+\varepsilon} |Z_t^{n,s} - Z_{t-\delta}^{n,s}|^2 \, ds \right] dt \right\} + \delta C$$

and consequently, applying Fubini's theorem and Fatou's lemma, and using (4.18), we obtain that

$$\varlimsup_{\delta \to 0} \sup_n E\left\{\int_0^T |Z_t^n - Z_{t-\delta}^n|^2 \, dt\right\}$$
$$\leq C \varlimsup_{\delta \to 0} \sup_n E\left\{\int_0^T |Z_t^{n,\infty} - Z_{t-\delta}^{n,\infty}|^2 \, dt\right\}$$
$$+ C \int_\varepsilon^T \varlimsup_{\delta \to 0} \sup_n E\left\{\int_0^{s-\varepsilon} |Z_t^{n,s} - Z_{t-\delta}^{n,s}|^2 \, dt\right\} ds + C\varepsilon = C\varepsilon.$$

Since $\varepsilon$ is arbitrary, we get

(4.19) $$\varlimsup_{\delta \to 0} \sup_n E\left\{\int_0^T |Z_t^n - Z_{t-\delta}^n|^2 \, dt\right\} = 0.$$

Now, note that

$$E \int_0^T |Z_t^n - Z_t^{n,\delta}|^2 \, dt \leq \frac{1}{\delta} E \int_0^T \int_{t-\delta}^t |Z_s^n - Z_t^n|^2 \, ds \, dt$$
$$= E \int_0^1 \int_0^T |Z_{t-\delta r}^n - Z_t^n|^2 \, dt \, dr,$$

we see that (4.19) implies (3.1), and the existence of the weak solution to (4.1) follows from Theorem 3.1.

Finally, from the proof of Theorem 3.1, we see that $(X^n, Y^n) \to (X, Y)$ in distribution. Thus, applying the Skorohod representation theorem again if necessary, we can assume without loss of generality that the convergence is $\mathbb{P}$-a.s. Note that $Y_t^n = u^n(t, X_t^n)$ and $u^n \to u$ uniformly, it follows that $Y_t = u(t, X_t)$, for all $t \in [0, T]$, $\mathbb{P}$-a.s. □

Finally, we prove a regularity of the weak solution above, which will be useful in the next section.

COROLLARY 4.6. *Assume* (H1), (H2) *and* $T_0 < T$. *Let* $(\Omega, \mathcal{F}, \mathbb{P}; \mathbf{F}, W, X, Y, Z)$ *be the weak solution constructed in Theorem* 4.5. *Then for any* $t \leq T_0$ *and any* $0 < \delta \leq \frac{T-T_0}{2}$, *it holds that*

(4.20) $$E_t^{\mathbb{P}}\left\{|Y_{t+\delta} - Y_t|^2 + \int_t^{t+\delta} |Z_s|^2 \, ds\right\} \leq \frac{C}{(T-T_0)^\alpha}\delta^\alpha, \qquad \mathbb{P}\text{-}a.s.,$$

*where* $E_t^{\mathbb{P}} \triangleq E^{\mathbb{P}}\{\cdot | \mathcal{F}_t\}$ *and* $\alpha$ *is the constant in Lemma* 4.1.

PROOF. Again let us denote $C > 0$ to be a generic constant which is allowed to vary from line to line. Applying Lemma 4.1 we have

$$|Y_{t+\delta} - Y_t| = |u(t+\delta, X_{t_\delta}) - u(t, X_t)|$$
$$\leq \frac{C}{(T-T_0)^{\alpha/2}}[\delta^{\alpha/2} + |X_{t+\delta} - X_t|^\alpha].$$

Therefore, using the boundedness of $\sigma$, we deduce

$$E_t^{\mathbb{P}}\{|Y_{t+\delta} - Y_t|^2\}$$

(4.21)
$$\leq \frac{C}{(T-T_0)^\alpha}\left[\delta^\alpha + E_t^{\mathbb{P}}\left\{\left|\int_t^{t+\delta} \sigma(s, X_s, Y_s)\, dW_s\right|^{2\alpha}\right\}\right]$$

$$\leq \frac{C}{(T-T_0)^\alpha}\left[\delta^\alpha + E_t^{\mathbb{P}}\left\{\left|\int_t^{t+\delta} |\sigma(s, X_s, Y_s)|^2\, ds\right|^{\alpha}\right\}\right]$$

$$\leq \frac{C}{(T-T_0)^\alpha}\delta^\alpha.$$

Moreover, note that

$$E_t^{\mathbb{P}}\left\{\int_t^{t+\delta}|Z_s|^2\, ds\right\} = E_t^{\mathbb{P}}\left\{\left|\int_t^{t+\delta} Z_s\, dW_s\right|^2\right\}$$

$$= E_t^{\mathbb{P}}\left\{\left|Y_{t+\delta} - Y_t + \int_t^{t+\delta} h(s, X_s, Y_s, Z_s)\, ds\right|^2\right\}$$

$$\leq 2E_t^{\mathbb{P}}\{|Y_{t+\delta} - Y_t|^2\} + C\delta^2.$$

This, together with (4.21), proves the (4.20). □

*Extension to the cases of "discrete functionals."* The result of Theorem 4.5 can be easily extended to FBSDEs whose coefficients are discrete functionals. For example, let $\pi: 0 = t_0 < t_1 < \cdots < t_N = T$ be a given partition, and denote $(X)_t \stackrel{\triangle}{=} (X_{t_1 \wedge t}, \ldots, X_{t_N \wedge t})$ for $t \in [0, T]$. Consider the following FBSDE:

(4.22)
$$\begin{cases} X_t = x + \int_0^t \sigma(s, (X)_s, Y_s)\, dW_s; \\ Y_t = g((X)_T) + \int_t^T h(s, (X)_s, Y_s, Z_s)\, ds, \qquad t \in [0, T], \\ \quad - \int_t^T Z_s\, dW_s; \end{cases}$$

where $h \in \mathbb{C}([0, T] \times \mathbb{R}^N \times \mathbb{R} \times \mathbb{R}; \mathbb{R})$, $g \in \mathbb{C}(\mathbb{R}^N; \mathbb{R})$, and both are uniformly bounded.

THEOREM 4.7. *Assume* (H1) *and* (H2). *Then* (4.22) *admits a weak solution.*

PROOF. Let $(\sigma_n, h_n, g_n)$ be smooth mollifiers of $(\sigma, h, g)$. We define functions $u_k^n$ backwardly as follows. First, considering $(x_1, \ldots, x_N)$ as parameters, we define

$$u_{N+1}^n(x_1, \ldots, x_N; T, x) \stackrel{\triangle}{=} g_n(x_1, \ldots, x_N).$$

For $k = N, \ldots, 1$, given $(x_1, \ldots, x_{k-1})$ as parameters, let $u_k^n(x_1, \ldots, x_{k-1}; t, x)$ be the classical solution to the following PDE over $t \in [t_{k-1}, t_k]$:

$$\begin{cases} \partial_t u_k^n + \frac{1}{2}\sigma_n^2(t, x_1, \ldots, x_{k-1}, \underbrace{x, \ldots, x}_{N-k}, u_k^n)\partial_{xx} u_k^n \\ \qquad + h_n(t, x_1, \ldots, x_{k-1}, \underbrace{x, \ldots, x}_{N-k}, u_k^n, \sigma_n \partial_x u_k^n) = 0; \\ u_k^n(x_1, \ldots, x_{k-1}; t_k, x) = u_{k+1}^n(x_1, \ldots, x_{k-1}, x; t_k, x). \end{cases}$$

Now, we construct the solutions $(X^n, Y^n, Z^n)$ recursively as follows. Define $X_0^n \triangleq x$. For $k = 1, \ldots, N$ and $t \in (t_{k-1}, t_k]$, we define

$$X_t^n = X_{t_{k-1}}^n + \int_{t_{i-1}}^t \sigma(s, (X^n)_{t_{k-1}}, X_s^n, u_n((X^n)_{t_{k-1}}; s, X_s^n))\, dW_s.$$

For $t \in [t_{k-1}, t_k)$, let

$$Y_t^n \triangleq u_k^n((X^n)_{t_{k-1}}; t, X_t^n); \qquad Z_t^n \triangleq \partial_x u_k^n((X^n)_{t_{k-1}}; t, X_t^n)\sigma_n(t, (X^n)_t, Y_t^n).$$

Then $(X^n, Y^n, Z^n)$ is a strong solution to FBSDE (4.22) with coefficients $(\sigma_n, h_n, g_n)$.

For $t \in [t_{N-1}, t_N]$, note that $(X^n, Y^n, Z^n)$ satisfies

$$\begin{cases} X_t^n = X_{t_{N-1}}^n + \int_{t_{N-1}}^t \sigma_n(s, (X^n)_{t_{N-1}}, X_s^n, Y_s^n)\, dW_s; \\ Y_t^n = u_n(X_{t_{N-1}}^n; T, X_{t_N}^n) \\ \qquad + \int_t^{t_N} h_n(s, (X^n)_{t_{N-1}}, X_s^n, Y_s^n, Z_s^n)\, ds - \int_t^{t_N} Z_s^n\, dW_s. \end{cases}$$

Since $u_N^n(x_1, \ldots, x_{N-1}; t_N, x)$ is uniformly continuous on $x$, following the same arguments as in previous subsection, we get

$$\lim_{\delta \to 0} \sup_n E\left\{ \int_{t_{N-1}}^{t_N} |Z_t^n - Z_t^{n,\delta}|^2\, dt \right\} = 0.$$

Moreover, by stability of the PDEs, we know $u_{N_1}^n(x_1, \ldots, x_{N-2}; t_{N-1}, x)$ is uniformly continuous on $x$, then we may prove

$$\lim_{\delta \to 0} \sup_n E\left\{ \int_{t_{N-2}}^{t_{N-1}} |Z_t^n - Z_t^{n,\delta}|^2\, dt \right\} = 0.$$

Repeat the arguments, we get

$$\lim_{\delta \to 0} \sup_n E\left\{ \int_{t_{k-1}}^{t_k} |Z_t^n - Z_t^{n,\delta}|^2\, dt \right\} = 0 \qquad \forall k.$$

Thus,

$$\lim_{\delta \to 0} \sup_n E\left\{ \int_0^T |Z_t^n - Z_t^{n,\delta}|^2\, dt \right\} = 0.$$

Now the result follows from Theorem 3.1 immediately. □

REMARK 4.8. We should point out that a decoupled version of (4.22) was studied in Hu and Ma [12], in which the existence of strong solution was proved under the assumption that $\sigma$ is Lipschitz. However, it is by no means clear if the method there can be extended to the current case.

**5. Uniqueness.** We now turn our attention to the key issue of the paper: the uniqueness of the solution to FBMP. Again, we shall consider only the special case (4.1) and assume all processes are one-dimensional. We shall assume throughout this section that (H1) and (H2) hold.

Recall from Section 2 the canonical space $\Omega \stackrel{\triangle}{=} \mathbb{C}([0,T];\mathbb{R}) \times \mathbb{C}([0,T];\mathbb{R})$. Let $\mathbf{F} = \{\mathcal{F}_t\}_{t\geq 0}$ denote the filtration generated by the canonical processes, which we shall denote by $(\mathbf{x}, \mathbf{y})$. In light of Remark 2.6(ii), from now on we call $(\mathbb{P}, Z)$ a solution to the $\text{FBMP}_{x,T}(0, \sigma, h, g)$. For simplicity, in what follows, we do not distinguish the term "solution to the FBMP" from "weak solution," and we often simply write "FBMP (4.1)" instead of "$\text{FBMP}_{x,T}(0, \sigma, h, g)$" when the context is clear.

We first give the definition of the uniqueness for FBMP.

DEFINITION 5.1. We say that the forward–backward martingale problem FBMP (4.1) has unique solution whenever $(\mathbb{P}^i, Z^i)$, $i = 1, 2$, are two solutions to the FBMP such that $\mathbb{P}^i(\mathbf{x}_0 = x) = 1$, $i = 1, 2$, then the processes $(\mathbf{x}, \mathbf{y}, Z^1)$ and $(\mathbf{x}, \mathbf{y}, Z^2)$ have the same finite dimensional distributions, under $\mathbb{P}^1$ and $\mathbb{P}^2$, respectively. In particular, this means that $\mathbb{P}^1 = \mathbb{P}^2$.

By the proof of Theorem 2.5, it is obvious that the uniqueness of solution to FBMP (4.1) is equivalent to the uniqueness in law of weak solution to FBSDE (4.1).

By Theorem 4.5, we know that there exists at least one solution to the FBMP (4.1). We denote this solution by $(\mathbb{P}^0, Z^0)$. We note that this special weak solution has the following feature:

$$\mathbf{y}_t = u(t, \mathbf{x}_t) \qquad \forall t \in [0, T], \ \mathbb{P}^0\text{-a.s.} \tag{5.1}$$

where $u$ is a viscosity solution to PDE (4.2) satisfying (4.3). Clearly, to prove the uniqueness of solution to FBMP (4.1), it suffices to show that any solution to it will be identical in law to $(\mathbb{P}^0, Z^0)$.

To begin with, we recall that for any given probability measure $\mathbb{P} \in \mathcal{P}(\Omega)$ and any $t < T$, there exists a *regular conditional probability distribution* (r.c.p.d. for short) of $\mathbb{P}$ given $\mathcal{F}_t$, denoted by $\mathbb{P}_t^\omega$, $\omega \in \Omega$, in the sequel (see, e.g., [26]). Furthermore, we can choose a version of $\mathbb{P}_t^\omega$ so that $\mathbb{P}_t^\omega \in \mathcal{P}(\Omega)$ for all $\omega \in \Omega$. In what follows, we will always take such a version without further specification.

We now introduce an auxiliary notion that will play an important role in our discussion for uniqueness. Let $k = k(t, \delta, \eta)$ be a (deterministic) function defined on $[0, T) \times (0, T) \times (0, 1)$ satisfying the following properties.

(5.2) $$\begin{cases} k(t_1, \delta_1, \eta) \leq k(t_2, \delta_2, \eta), & \forall t_1 \leq t_2, \delta_1 \leq \delta_2; \\ \lim_{\delta \to 0} k(t, \delta, \eta) = 0, & \forall (t, \eta); \\ k(t, \delta, \eta) \geq \dfrac{C\delta^\alpha}{(T-t-\delta)^\alpha \eta^2}, & \forall t + \delta < T, \end{cases}$$

where $C, \alpha > 0$ are the constants same as those in Corollary 4.6.

DEFINITION 5.2. We say that a weak solution $(\mathbb{P}, Z)$ is a "$k$-weak solution" at $(t, x, y)$ if the following hold:

(i) $W_s \stackrel{\triangle}{=} \int_t^s \sigma^{-1}(r, \mathbf{x}_r, \mathbf{y}_r) \, d\mathbf{x}_r$ is a $\mathbb{P}$-Brownian motion for $s \geq t$;
(ii) $\mathbb{P}\{\mathbf{x}_t = x, \ \mathbf{y}_t = y\} = 1$;
(iii) $\mathbf{y}_s = y - \int_t^s h(r, \mathbf{x}_r, \mathbf{y}_r, Z_r) \, dr + \int_t^s Z_r \, dW_r$, $s \in [t, T]$, $\mathbb{P}$-a.s.;
(iv) $\mathbb{P}\{\mathbf{y}_T = g(\mathbf{x}_T)\} = 1$;
(v) for any $s \in [t, T)$, and $\delta > 0$, $\eta > 0$, $\mathbb{P}_s^\omega\{|\mathbf{y}_s - \mathbf{y}_{(s+\delta)\wedge T}| \geq \eta\} \leq k(s, \delta, \eta)$, $\mathbb{P}$-a.s. $\omega \in \Omega$.

REMARK 5.3. (i) Note that in light of Corollary 4.6 one can easily show that

(5.3) $$\mathbb{P}_t^{0,\omega}\{|\mathbf{y}_t - \mathbf{y}_{t+\delta}| \geq \eta\} \leq \frac{C\delta^\alpha}{(T-t-\delta)^\alpha \eta^2}, \qquad \mathbb{P}^0\text{-a.s. } \omega \in \Omega.$$

Therefore, for a given function $k$ satisfying (5.2), any weak solution constructed via solution to the PDE with initial time $t$ will be a $k$-weak solution at $(t, \mathbf{x}_t, u(t, \mathbf{x}_t))$. In particular, the solution $(\mathbb{P}^0, Z^0)$ is a $k$-weak solution at $(0, x, u(0, x))$.

(ii) Note that if there exists a function $\tilde{k}$ such that $\tilde{k}(t_1, \delta_1) \leq \tilde{k}(t_2, \delta_2)$ for $t_1 \leq t_2$, $\delta_1 \leq \delta_2$ and $\lim_{\delta \to 0} \tilde{k}(s, \delta) = 0$ for $t < T$, and that for any $s \in [t, T)$, and $\delta > 0$,

(5.4) $$E^{\mathbb{P}_s^\omega}\left\{\int_s^{(s+\delta)\wedge T} |Z_r|^2 \, dr\right\} \leq \tilde{k}(s, \delta), \qquad \mathbb{P}\text{-a.s. } \omega \in \Omega,$$

then it is easy to check, by the Burkholder–Davis–Gundy and Hölder inequalities that condition (v) in Definition 5.2 holds for an appropriately chosen function $k$.

(iii) Since $h$ and $g$ are bounded, by conditions (iii) and (iv), one can easily show that $E^{\mathbb{P}}\{\int_t^T |Z_s|^2 \, ds\} < \infty$. $Z$ can be arbitrary over $[0, t)$, as long as $\mathbb{P}$ and $\mathbf{z}_s \stackrel{\triangle}{=} Z_s \sigma^{-1}(s, \mathbf{x}_s, \mathbf{y}_s)$ satisfy the conditions in Definition 2.3.

We shall prove that any $k$-weak solution is identical to $(\mathbb{P}^0, Z^0)$, by showing that a $k$-weak solution can exist only at $(t, x, u(t, x))$. To this end, let us denote

$$(5.5) \qquad \mathcal{O} \triangleq \{(t, x, y) : \text{there exists a } k\text{-weak solution at } (t, x, y)\},$$

and let $\bar{\mathcal{O}}$ denote the closure of $\mathcal{O}$ (we note that $\mathcal{O}$ is *not* necessarily Lebesgue measurable!). Clearly, we have $(t, x, u(t, x)) \in \mathcal{O}$ for any $(t, x) \in [0, T] \times \mathbb{R}$.

Now define two functions on $(t, x) \in [0, T] \times \mathbb{R}$:

$$(5.6) \quad \underline{u}(t, x) \triangleq \inf\{y : (t, x, y) \in \bar{\mathcal{O}}\}; \qquad \bar{u}(t, x) \triangleq \sup\{y : (t, x, y) \in \bar{\mathcal{O}}\}.$$

We claim that for some constant $C_0$,

$$(5.7) \quad -C_0 \leq \underline{u}(t, x) \leq u(t, x) \leq \bar{u}(t, x) \leq C_0; \qquad \underline{u}(T, x) = \bar{u}(T, x) = g(x).$$

First, for any $(t, x, y) \in \mathcal{O}$, let $(\mathbb{P}, Z)$ be a $k$-weak solution at $(t, x, y)$. Then

$$y = g(\mathbf{x}_T) - \int_t^T h(s, \mathbf{x}_s, \mathbf{y}_s, Z_s)\, ds + \int_t^T Z_s\, dW_s, \qquad \mathbb{P}\text{-a.s.}$$

Thus,

$$(5.8) \qquad y = E^{\mathbb{P}} y = E^{\mathbb{P}} \Big\{ g(\mathbf{x}_T) - \int_t^T h(s, \mathbf{x}_s, \mathbf{y}_s, Z_s)\, ds \Big\}.$$

Since both $g$ and $h$ are bounded, there exists some $C_0 > 0$ such that

$$|y|^2 \leq C E^{\mathbb{P}} \Big\{ |g(\mathbf{x}_T)|^2 + \int_t^T |h(s, \mathbf{x}_s, \mathbf{y}_s, Z_s)|^2\, ds \Big\} \leq C_0^2.$$

Second, for any $(T, x, y) \in \bar{\mathcal{O}}$, assume $(t_n, x_n, y_n) \in \mathcal{O}$ and $(t_n, x_n, y_n) \to (T, x, y)$. Let $(\mathbb{P}^n, Z^n)$ be a $k$-weak solution at $(t_n, x_n, y_n)$ and $W^n$ be the corresponding $\mathbb{P}^n$-Brownian motion. Then

$$\begin{cases} \mathbf{x}_T = x_n + \displaystyle\int_{t_n}^T \sigma(s, \mathbf{x}_s, \mathbf{y}_s)\, dW_s^n; \\ y_n = g(\mathbf{x}_T) - \displaystyle\int_{t_n}^T h(s, \mathbf{x}_s, \mathbf{y}_s, Z_s)\, ds + \displaystyle\int_{t_n}^T Z_s\, dW_s^n; \end{cases} \quad \mathbb{P}^n\text{-a.s.}$$

Thus, applying the similar argument as (5.8) from the FBSDE above we have

$$|y_n - g(x_n)|^2$$
$$\leq 2 E^{\mathbb{P}^n} \Big\{ |g(\mathbf{x}_T) - g(x_n)|^2 + \Big| \int_{t_n}^T h(s, \mathbf{x}_s, \mathbf{y}_s, Z_s)\, ds \Big|^2 \Big\}$$
$$\leq C E^{\mathbb{P}^n} \Big\{ \Big| g\Big(x_n + \int_{t_n}^T \sigma(s, \mathbf{x}_s, \mathbf{y}_s)\, dW_s^n\Big) - g(x_n) \Big|^2 \Big\} + C|T - t_n|^2.$$

Now note that by (H1) $g$ is bounded and uniformly continuous, a standard argument using Chebyshev's inequality and the boundedness of $\sigma$, one shows easily that $\lim_{n\to\infty} |y_n - g(x_n)| = 0$. To wit, $y = g(x)$.

Moreover, since $\bar{\mathcal{O}}$ is a closed set, we have $(t, x, \underline{u}(t,x)) \in \bar{\mathcal{O}}$ and $(t, x, \bar{u}(t,x)) \in \bar{\mathcal{O}}$.

LEMMA 5.4. *$\underline{u}$ is lower semi-continuous and $\bar{u}$ is upper semi-continuous.*

PROOF. We need only check for $\underline{u}$. The argument for $\bar{u}$ is symmetric. Assume $(t_n, x_n) \to (t_0, x_0)$ and $\underline{u}(t_n, x_n) \to y_0$. Since $(t_n, x_n, \underline{u}(t_n, x_n)) \in \bar{\mathcal{O}}$ for each $n$ and $\bar{\mathcal{O}}$ is closed, we see that $(t_0, x_0, y_0) \in \bar{\mathcal{O}}$, hence $\underline{u}(t_0, x_0) \leq y_0$. □

Our main result of this section is the following theorem.

THEOREM 5.5. *Assume* (H1), (H2) *and* (H3). *Then, $\underline{u}$ and $\bar{u}$ are viscosity supersolution and subsolution, respectively, of the quasilinear PDE*

$$(5.9) \quad \begin{cases} u_t + \frac{1}{2}\sigma^2(t,x,u)u_{xx} + h(t,x,u,u_x\sigma) = 0; \\ u(T,x) = g(x). \end{cases}$$

PROOF. Again, we check for $\underline{u}$ only. For any $(t_0, x_0) \in [0,T) \times \mathbb{R}$, let $\varphi \in C^{1,2}([0,T] \times \mathbb{R})$ be such that $y_0 \stackrel{\triangle}{=} \underline{u}(t_0, x_0) = \varphi(t_0, x_0)$ and $\underline{u}(t,x) \geq \varphi(t,x)$, for all $(t,x) \in [0,T] \times \mathbb{R}$. We shall prove that

$$(5.10) \quad [\mathcal{L}\varphi](t_0, x_0, \varphi(t_0, x_0)) \leq 0,$$

where

$$[\mathcal{L}\varphi](t,x,y) \stackrel{\triangle}{=} \varphi_t(t,x) + \tfrac{1}{2}\sigma^2(t,x,y)\varphi_{xx}(t,x)$$
$$+ h(t,x,y,\varphi_x(t,x)\sigma(t,x,y)).$$

To do this, we first note that $(t_0, x_0, y_0) = (t_0, x_0, \underline{u}(t_0, x_0)) \in \bar{\mathcal{O}}$, so for each $n$ there exists a $(t_n, x_n, y_n) \in \mathcal{O}$ such that

$$(5.11) \quad |t_n - t_0| + |x_n - x_0| + |y_n - y_0| \leq \frac{1}{n}.$$

Now suppose that $(\mathbb{P}^n, Z^n)$ is a $k$-weak solution at $(t_n, x_n, y_n)$ and let $W^n$ denote the corresponding $\mathbb{P}^n$-Brownian motion. For $t \in (t_n, T)$, it is readily seen that $(P_t^{n,\omega}, Z)$ is a $k$-weak solution at $(t, \mathbf{x}_t, \mathbf{y}_t)$, $\mathbb{P}^n$-a.s. $\omega \in \Omega$. In other words, we must have $(t, \mathbf{x}_t, \mathbf{y}_t) \in \mathcal{O}$, $\mathbb{P}^n$-a.s., and consequently $\mathbf{y}_t \geq \underline{u}(t, \mathbf{x}_t) \geq \varphi(t, \mathbf{x}_t)$, $\mathbb{P}^n$-a.s., $\forall t \geq t_n$.

Now let us denote

$$\Delta Y_t \stackrel{\triangle}{=} \varphi(t, \mathbf{x}_t) - \mathbf{y}_t; \qquad \Delta Z_t^n \stackrel{\triangle}{=} \varphi_x \sigma(t, \mathbf{x}_t, \mathbf{y}_t) - Z_t^n.$$

Also, for any $\varepsilon > 0$, let $h_\varepsilon$ be a mollifier of $h$ such that $\|h_\varepsilon - h\|_\infty \leq \varepsilon$ and $\|\partial_z h_\varepsilon\|_\infty \leq C_\varepsilon$, and denote

$$\alpha_t^{n,\varepsilon} \triangleq [h(t, \mathbf{x}_t, \mathbf{y}_t, Z_t^n) - h(t, \mathbf{x}_t, \mathbf{y}_t, \varphi_x \sigma(t, \mathbf{x}_t, \mathbf{y}_t))]$$
$$\quad - [h_\varepsilon(t, \mathbf{x}_t, \mathbf{y}_t, Z_t^n) - h_\varepsilon(t, \mathbf{x}_t, \mathbf{y}_t, \sigma(t, \mathbf{x}_t, \mathbf{y}_t)\varphi_x)];$$
$$\beta_t^{n,\varepsilon} \triangleq \int_0^1 \partial_z h_\varepsilon(t, \mathbf{x}_t, \mathbf{y}_t, Z_t^n + \theta \Delta Z_t^n)\, d\theta.$$

Then it holds that

(5.12) $$|\alpha_t^{n,\varepsilon}| \leq 2\varepsilon, \qquad |\beta_t^{n,\varepsilon}| \leq C_\varepsilon.$$

Furthermore, applying Itô's formula and using the definition of $\mathcal{L}\varphi$, $\alpha^{n,\varepsilon}$ and $\beta^{n,\varepsilon}$ we have

$$d\Delta Y_t = [\varphi_t + \tfrac{1}{2}\sigma^2(t, \mathbf{x}_t, \mathbf{y}_t)\varphi_{xx} + h(t, \mathbf{x}_t, \mathbf{y}_t, Z_t^n)]\, dt + \Delta Z_t^n\, dW_t^n$$
$$= \{[\mathcal{L}\varphi](t, \mathbf{x}_t, \mathbf{y}_t) + [h(t, \mathbf{x}_t, \mathbf{y}_t, Z_t^n) - h(t, \mathbf{x}_t, \mathbf{y}_t, \varphi_x \sigma(t, \mathbf{x}_t, \mathbf{y}_t))]\}\, dt$$
$$\quad + \Delta Z_t^n\, dW_t^n$$
$$= [\mathcal{L}\varphi](t, \mathbf{x}_t, \mathbf{y}_t)\, dt + \alpha_t^{n,\varepsilon}\, dt - \beta_t^{n,\varepsilon} \Delta Z_t^n\, dt + \Delta Z_t^n\, dW_t^n.$$

Now let us denote

$$\Gamma_t^{n,\varepsilon} \triangleq \exp\left\{\int_{t_n}^t \beta_s^{n,\varepsilon}\, dW_s^n - \tfrac{1}{2}\int_{t_n}^t |\beta_s^{n,\varepsilon}|^2\, ds\right\}, \qquad t \in [t_n, T].$$

One easily checks that by denoting $E^n \triangleq E^{\mathbb{P}^n}$,

(5.13)
$$\Gamma_{t_n}^{n,\varepsilon} = 1, \qquad \Gamma_t^{n,\varepsilon} > 0, \qquad E^n\{\Gamma_t^{n,\varepsilon}\} = 1 \quad \text{and}$$
$$E^n\{|\Gamma_t^{n,\varepsilon}|^2\} \leq C_\varepsilon \qquad \forall t \geq t_n.$$

Moreover, applying Itô's formula again we have, for $t \in [t_n, T]$,

(5.14)
$$d(\Gamma_t^{n,\varepsilon} \Delta Y_t) = \Gamma_t^{n,\varepsilon}[\mathcal{L}\varphi]\, dt + \Gamma_t^{n,\varepsilon} \alpha_t^{n,\varepsilon}\, dt$$
$$\quad + \Gamma_t^{n,\varepsilon}[\Delta Z_t^n - \beta_t^{n,\varepsilon}\Delta Y_t]\, dW_t^n.$$

Now, for any $\delta > \frac{1}{n}$, choose $t = t_0 + \delta > t_n$ [see (5.11)], we deduce from (5.14) that

$$0 \geq E^n\{\Gamma_{t_0+\delta}^{n,\varepsilon} \Delta Y_{t_0+\delta}\} = E^n\left\{\Delta Y_{t_n} + \int_{t_n}^{t_0+\delta} \Gamma_t^{n,\varepsilon}\{[\mathcal{L}\varphi](t, \mathbf{x}_t, \mathbf{y}_t) + \alpha_t^{n,\varepsilon}\}\, dt\right\}.$$

Therefore, using (5.12) and (5.13), we get

$$E^n\left\{\int_{t_n}^{t_0+\delta} \Gamma_t^{n,\varepsilon}[\mathcal{L}\varphi](t, \mathbf{x}_t, \mathbf{y}_t)\, dt\right\}$$

$$\leq -E^n \left\{ \Delta Y_{t_n} + \int_{t_n}^{t_0+\delta} \Gamma_t^{n,\varepsilon} \alpha_t^{n,\varepsilon} \, dt \right\}$$

(5.15)
$$\leq E^n \left\{ |y_n - y_0| + |\varphi(t_0, x_0) - \varphi(t_n, x_n)| + \int_{t_n}^{t_0+\delta} \Gamma_t^{n,\varepsilon} |\alpha_t^{n,\varepsilon}| \, dt \right\}$$

$$\leq CE^n \left\{ \frac{1}{n} + \varepsilon \int_{t_n}^{t_0+\delta} \Gamma_t^{n,\varepsilon} \, dt \right\} \leq C\left[\varepsilon + \frac{1}{n\delta - 1}\right](t_0 + \delta - t_n),$$

where $C$ may depend on $\varphi$. Recall (5.7). To prove (5.10), without loss of generality, we may assume that $\varphi(t,x) = -C_0 - 1$ for $x$ outside of a compact set. Then $\varphi, \varphi_t, \varphi_x$ and $\varphi_{xx}$ are all uniformly continuous. By (H1) and (H3), $\mathcal{L}\varphi$ is uniformly continuous in $(t,x,y)$. Let $w_\varphi$ denote the modulus of continuity of $\mathcal{L}\varphi$, and write

$$\Delta_n[\mathcal{L}\varphi](t, \mathbf{x}_t, \mathbf{y}_t) = \mathcal{L}\varphi(t, \mathbf{x}_t, \mathbf{y}_t) - \mathcal{L}\varphi(t_n, x_n, y_n).$$

We see that (5.15) yields

$$\mathcal{L}\varphi(t_0, x_0, y_0)$$
$$\leq \mathcal{L}\varphi(t_n, x_n, y_n) + w_\varphi\left(\frac{1}{n}\right)$$
$$= E^n \left\{ \frac{1}{t_0 + \delta - t_n} \int_{t_n}^{t_0+\delta} \Gamma_t^{n,\varepsilon} \mathcal{L}\varphi(t_n, x_n, y_n) \, dt \right\} + w_\varphi\left(\frac{1}{n}\right)$$

(5.16) $$= E^n \left\{ \frac{1}{t_0 + \delta - t_n} \int_{t_n}^{t_0+\delta} \Gamma_t^{n,\varepsilon} \{[\mathcal{L}\varphi](t, \mathbf{x}_t, \mathbf{y}_t) - \Delta_n[\mathcal{L}\varphi](t, \mathbf{x}_t, \mathbf{y}_t)\} \, dt \right\}$$
$$+ w_\varphi\left(\frac{1}{n}\right)$$
$$\leq C\varepsilon + \frac{C}{n\delta - 1} + w_\varphi\left(\frac{1}{n}\right)$$
$$+ \frac{1}{t_0 + \delta - t_n} E^n \left\{ \int_{t_n}^{t_0+\delta} |\Gamma_t^{n,\varepsilon} \Delta_n[\mathcal{L}\varphi](t, \mathbf{x}_t, \mathbf{y}_t)| \, dt \right\}.$$

To estimate the last term on the right-hand side above, we first apply the Cauchy–Schwarz inequality and the estimate (5.13) to get

$$E^n \int_{t_n}^{t_0+\delta} |\Gamma_t^{n,\varepsilon} \Delta_n[\mathcal{L}\varphi]| \, dt$$

(5.17) $$\leq \left\{ E^n \int_{t_n}^{t_0+\delta} |\Gamma_t^{n,\varepsilon}|^2 \, dt \right\}^{1/2} \left\{ E^n \int_{t_n}^{t_0+\delta} |\Delta_n[\mathcal{L}\varphi]|^2 \, dt \right\}^{1/2}$$
$$\leq C_\varepsilon \left\{ \sup_{t_n \leq t \leq t_0+\delta} E^n\{|\Delta_n[\mathcal{L}\varphi](t, \mathbf{x}_t, \mathbf{y}_t)|^2\} \right\}^{1/2} (t_0 + \delta - t_n).$$

Note that, for any $\eta > 0$ and $t \in [t_n, t_0 + \delta]$, we apply the Chebyshev inequality to get

$$
\begin{aligned}
E^n\{|\Delta_n[\mathcal{L}\varphi](t, \mathbf{x}_t, \mathbf{y}_t)|^2\} \\
&\leq Cw_\varphi^2(t_0 + \delta - t_n) + Cw_\varphi^2(\eta) \\
&\quad + CP^n(|\mathbf{x}_t - x_n| \geq \eta) + CP^n(|\mathbf{y}_t - y_n| \geq \eta) \\
&\leq C\bigg[w_\varphi^2(t_0 + \delta - t_n) + w_\varphi^2(\eta) \\
&\quad + \frac{1}{\eta^2}E^n\bigg\{\int_{t_n}^t |\sigma(s, \mathbf{x}_s, \mathbf{y}_s)|^2\,ds\bigg\} + P^n(|\mathbf{y}_t - \mathbf{y}_{t_n}| \geq \eta)\bigg] \\
&\leq C\bigg[w_\varphi^2(t_0 + \delta - t_n) + w_\varphi^2(\eta) \\
&\quad + \frac{1}{\eta^2}[t_0 + \delta - t_n] + k(t_n, t_0 + \delta - t_n, \eta)\bigg],
\end{aligned}
\tag{5.18}
$$

thanks to Definition 5.2(v) and (5.2). Combining (5.17) and (5.18), we see that (5.16) now becomes

$$
\begin{aligned}
\mathcal{L}\varphi(t_0, x_0, y_0) \\
&\leq C\varepsilon + \frac{C}{n\delta - 1} + w_\varphi\bigg(\frac{1}{n}\bigg) \\
&\quad + C_\varepsilon\bigg[w_\varphi(t_0 + \delta - t_n) + w_\varphi(\eta) \\
&\quad + \frac{1}{\eta}[t_0 + \delta - t_n]^{1/2} + k^{1/2}(t_n, t_0 + \delta - t_n, \eta)\bigg].
\end{aligned}
\tag{5.19}
$$

Now fix $\varepsilon$ and $\eta$, choose $\delta = \frac{1}{\sqrt{n}}$, and let $n \to \infty$. By (5.2), we get

$$\mathcal{L}\varphi(t_0, x_0, y_0) \leq C\varepsilon + C_\varepsilon w_\varphi(\eta).$$

Finally, letting $\eta \to 0$ and then $\varepsilon \to 0$, we obtain (5.10). That is, $\underline{u}$ is a viscosity supersolution, proving the theorem. □

A direct consequence of Theorem 5.5 is the following uniqueness result.

THEOREM 5.6. *Assume* (H1), (H2) *and* (H3), *and that the comparison theorem holds for bounded viscosity solutions to the PDE (4.2). Then FB-SDE (4.1) admits a unique weak solution* $(\mathbb{P}, Z)$ *satisfying* $\mathbb{P}_t^\omega\{|\mathbf{y}_t - \mathbf{y}_{(t+\delta)\wedge T}| \geq \eta\} \leq k(t, \delta, \eta), \mathbb{P}$-*a.s.* $\omega \in \Omega$, *for any* $(t, \delta, \eta)$.

PROOF. It suffices to show that $(\mathbb{P}, Z)$ is identical to the "canonical" weak solution $(\mathbb{P}^0, Z^0)$ constructed in Section 4, in the sense of Definition 5.1. We shall assume without loss of generality that $\mathbb{P}(\mathbf{y}_0 = y) = 1$ for some $y$ (otherwise, we apply the usual arguments by considering the conditional probabilities $\mathbb{P}^y\{\cdot\} \triangleq \mathbb{P}\{\cdot|\mathbf{y}_0 = y\}$, for $\mathbb{P}$-a.e. $y \in \mathbb{R}$, and the result will be the same). Then $(\mathbb{P}^0, Z^0)$ and $(\mathbb{P}, Z)$ are $k$-weak solutions at $(0, x, u(0, x))$ and $(0, x, y)$, respectively. Since by (5.7) and Theorem 5.5, we know that $\bar{u}$ is a bounded subsolution and $\underline{u}$ is a bounded supersolution to (4.2), by our assumptions we must have $\bar{u} \leq \underline{u}$, thanks to the comparison theorem. Thus, we must have $\underline{u} = u = \bar{u}$. On the other hand, following the arguments in Theorem 5.5, one shows that $(t, \mathbf{x}_t, \mathbf{y}_t) \in \bar{\mathcal{O}}$, $\mathbb{P}$-a.s., for any $t$. Therefore, it holds that $\underline{u}(t, \mathbf{x}_t) \leq \mathbf{y}_t \leq \bar{u}(t, \mathbf{x}_t)$. Thus, $\mathbf{y}_t = u(t, \mathbf{x}_t)$, $\mathbb{P}$-a.s., for all $t \in [0, T]$. Finally, since $\mathbf{x}, \mathbf{y}, u$ are continuous, we get $\mathbf{y}_t = u(t, \mathbf{x}_t)$, for all $t \in [0, T]$, $\mathbb{P}$-a.s.

Now define $dW \triangleq \sigma^{-1}(t, \mathbf{x}_t, u(t, \mathbf{x}_t)) \, d\mathbf{x}_t$, we see that $W$ is a $\mathbb{P}$-Brownian motion, and $(W, \mathbf{x})$ is a weak solution to a forward SDE. Since under our assumptions the uniqueness in law holds for this forward SDE, noting the relation $\mathbf{y}_t = u(t, \mathbf{x}_t)$ it is easily seen that $\mathbb{P} \circ (W, \mathbf{x}, \mathbf{y})^{-1} = \mathbb{P}^0 \circ (W, \mathbf{x}, \mathbf{y})^{-1}$. In particular, since $(\mathbf{x}, \mathbf{y})$ is the canonical process, we have $\mathbb{P} = \mathbb{P}^0$. Consequently, the processes $Z$ and $Z^0$, being the integrands of the quadratic variation processes $[\mathbf{y}, W]$ under $\mathbb{P}$ and $\mathbb{P}^0$, respectively, must be identical in law as well. In other words, the weak solutions $(\mathbb{P}, Z)$ and $(\mathbb{P}^0, Z^0)$ are identical by Definition 5.1, proving the theorem. □

REMARK 5.7. The assumption that the comparison theorem holds for the viscosity solution to the PDE (4.2) actually imply that the coefficients $\sigma$, $g$ and $h$ must satisfy certain conditions. We refer to the ubiquitous reference [5] for general theory of viscosity solutions. We should note that in general the sufficient conditions for comparison should be checked case by case, and we feel that it is more convenient to assume comparison theorem directly in Theorem 5.6. We shall present some simple cases in the concluding discussion below to make our point clearer.

We shall conclude our discussion on uniqueness by presenting some sufficient conditions under which the comparison theorem holds. We note that these cases are consequences of [5], Theorem 8.2, and the discussion in Section 5.D there, adjusted to the current situations. In light of the PDE (4.2), let us denote

$$F(t, x, y, p, A) \triangleq \tfrac{1}{2}\sigma^2(t, x, y)A + h(t, x, y, p\sigma(t, x, y)).$$

Suppose that:

(i) $F$ is decreasing in $y$;

(ii) there exists a continuous function $w:[0,\infty] \to [0,\infty]$, such that $w(0) = 0$ and for all $t$, $x_1$, $x_2$, $y$, $p$, $A$, $B$, and $\Delta x \stackrel{\triangle}{=} x_1 - x_2$; and all $\alpha > 0$ satisfying

$$-3\alpha \begin{bmatrix} I & 0 \\ 0 & I \end{bmatrix} \le \begin{bmatrix} A & 0 \\ 0 & -B \end{bmatrix} \le 3\alpha \begin{bmatrix} I & -I \\ -I & I \end{bmatrix},$$

it holds that

(5.20) $\quad |F(t,x_1,y,\alpha \Delta x, A) - F(t,x_2,y,\alpha \Delta x, B)| \le w(\alpha |\Delta x|^2 + |\Delta x|).$

Then it is known (cf. [5]) that the comparison theorem holds for viscosity subsolutions and supersolutions to PDE (4.2) that are of at most linear growth.

Now we assume that (H1)–(H3) hold. Since $\sigma$ is bounded, one sufficient condition for (ii) is $\sigma = \sigma(t,y)$ and $h$ is uniformly continuous in $x$ and $z$. Furthermore, in the following two examples, the condition (i) is satisfied as well, and consequently the comparison theorem holds.

EXAMPLE 1. $\sigma = \sigma(t)$ and $h$ is decreasing in $y$.

EXAMPLE 2. $\sigma = \sigma(t,y)$, $\sigma, h$ are uniformly Hölder continuous in $(x,y,z)$ and $g \in \mathbb{C}^2$ with bounded derivatives. Then (4.2) has a classical solution $u \in \mathbb{C}^{1,2}$. Assume $|u_x| \le C_1, |u_{xx}| \le C_2$ for some constants $C_1, C_2$. Assume further that, for any $t,x$ and $y_1 < y_2$,

$$\inf_{|p| \le C_1} [h(t,x,y_1,p\sigma(t,y_1)) - h(t,x,y_2,p\sigma(t,y_2))] \ge C_2 |\sigma(t,y_1) - \sigma(t,y_2)|.$$

Then one can check that the comparison theorem holds for all viscosity subsolutions and supersolutions that have bounded first and second derivatives in $x$.

**Acknowledgment.** The authors are very grateful to the anonymous referee for his/her extremely careful reading of the original manuscript and many helpful suggestions, which made this paper much better in both of its mathematical quality and its presentation.

J. MA
DEPARTMENT OF MATHEMATICS
UNIVERSITY OF SOUTHER CALIFORNIA
LOS ANGELES, CALIFORNIA 90089
USA
E-MAIL: jinma@usc.edu

J. ZHANG
DEPARTMENT OF MATHEMATICS
UNIVERSITY OF SOUTHERN CALIFORNIA
LOS ANGELES, CALIFORNIA 90089
USA
E-MAIL: jianfenz@usc.edu

Z. ZHENG
DEPARTMENT OF MATHEMATICAL SCIENCES
UNIVERSITY OF WISCONSIN–MILWAUKEE
BARCLAYS CAPITAL
745 7TH AVENUE
NEW YORK, NEW YORK 10019
USA
E-MAIL: ziyu.zheng@barclayscapital.com